\documentclass[11pt]{article}

\hsize \textwidth
\advance \hsize by -\marginparwidth
\evensidemargin \oddsidemargin
\usepackage{amssymb}
\advance\hoffset by 5mm

\usepackage{amssymb, amsmath,amsrefs,graphicx}

\newtheorem{df}{Definition}[section]
\newtheorem{corollary}[df]{Corollary}
\newtheorem{lemma}[df]{Lemma}
\newtheorem{prop}[df]{Proposition}
\newtheorem{thm}[df]{Theorem}

\newcommand{\farc}{\frac}

\makeatletter \@addtoreset{equation}{section}

\newcommand{\sgn}{\hbox{sgn}}
\newcommand{\Rm}{{\mathbb R}}
\newcommand{\Nm}{{\mathbb N}}

\def\vect#1{\mbox{\boldmath{$#1$}}}
\def\wG{\widehat{G}}
\def\wg{\widehat{g}}

\def\bfrho{\mbox{\boldmath$\rho$}}
\def\mR{\mathbb{R}}

\newcommand{\bes}{\begin{displaymath}}
\newcommand{\ees}{\end{displaymath}}
\newcommand{\be}{\begin{equation}}
\newcommand{\ee}{\end{equation}}
\newcommand{\ba}{\begin{eqnarray}}
\newcommand{\ea}{\end{eqnarray}}
\newcommand{\bas}{\begin{eqnarray*}}
\newcommand{\eas}{\end{eqnarray*}}
\newcommand{\@Bbb}[1]{\ensuremath{\mathbb #1}}

\newcommand{\B}{{\@Bbb B}}
\newcommand{\C}{{\@Bbb C}}

\newcommand{\F}{{\@Bbb F}}

\newcommand{\Q}{{\@Bbb Q}}
\newcommand{\bQ}{{\@Bbb Q}}
\newcommand{\N}{{\@Bbb N}}

\newcommand{\bbR}{{\@Bbb R}}
\newcommand{\W}{{\@Bbb W}}

\newcommand{\bbZ}{{\@Bbb Z}}
\newcommand{\bbT}{{\@Bbb T}}

\newcommand{\ltr}{\left\langle}
\newcommand{\rtr}{\right\rangle}

\newcommand{\lc}{\left(}
\newcommand{\rc}{\right)}

\newcommand{\lcu}{\left\{}
\newcommand{\rcu}{\right\}}

\newcommand{\eps}{\varepsilon}

\newcommand{\@s}[1]{\ensuremath{\mathcal #1}}
\newcommand{\cA}{\@s A}
\newcommand{\cB}{\@s B}
\newcommand{\cC}{\@s C}
\newcommand{\cD}{\@s D}
\newcommand{\cE}{\@s E}
\newcommand{\cF}{\@s F}
\newcommand{\cG}{\@s G}
\newcommand{\cH}{\@s H}
\newcommand{\cI}{\@s I}
\newcommand{\cJ}{\@s J}

\newcommand{\cK}{\@s K}
\newcommand{\cL}{\@s L}
\newcommand{\cN}{\@s N}
\newcommand{\cM}{\@s M}
\newcommand{\cO}{\@s O}
\newcommand{\cP}{\@s P}
\newcommand{\cR}{\@s R}
\newcommand{\cS}{\@s S}
\newcommand{\cT}{\@s T}
\newcommand{\cV}{\@s V}
\newcommand{\cW}{\@s W}
\newcommand{\cX}{\@s X}
\newcommand{\cY}{\@s Y}
\newcommand{\cZ}{\@s Z}

\newcommand{\@bm}[1]{\ensuremath{\mathbf #1}}
\newcommand{\bma}{\@bm a}
\newcommand{\bmb}{\@bm b}
\newcommand{\bmc}{\@bm c}
\newcommand{\bmd}{\@bm d}

\newcommand{\bme}{\@bm e}
\newcommand{\bmf}{\@bm f}
\newcommand{\bmg}{\@bm g}
\newcommand{\bmh}{\@bm h}
\newcommand{\bmi}{\@bm i}
\newcommand{\bmj}{\@bm j}
\newcommand{\bmk}{\@bm k}
\newcommand{\bml}{\@bm l}
\newcommand{\bmm}{\@bm m}
\newcommand{\bmn}{\@bm n}
\newcommand{\bmo}{\@bm o}
\newcommand{\bmp}{\@bm p}
\newcommand{\bmq}{\@bm q}
\newcommand{\bmr}{\@bm r}
\newcommand{\bms}{\@bm s}
\newcommand{\bmt}{\@bm t}
\newcommand{\bmu}{\@bm u}
\newcommand{\bmw}{\@bm w}
\newcommand{\bmv}{\@bm v}
\newcommand{\bmx}{\@bm x}
\newcommand{\bx}{\@bm x}

\newcommand{\bmy}{\@bm y}
\newcommand{\bz}{\@bm z}

\newcommand{\by}{\@bm y}

\newcommand{\bmzero}{\@bm 0}

\newcommand{\@g}[1]{\ensuremath{\mathfrak #1}}
\newcommand{\gA}{\@g A}
\newcommand{\gD}{\@g D}
\newcommand{\gJ}{\@g J}
\newcommand{\gF}{\@g F}
\newcommand{\gM}{\@g M}
\newcommand{\gR}{\mathbb R}

\newcommand{\argmin}{\mathop{\mbox{argmin}}}

\newcommand{\commentout}[1]{{}}

\makeatother

\title{A differential equations approach to $l_1$-minimization with applications to array imaging}
 \author{Miguel Moscoso\footnote{Escuela Polit\'ecnica Superior
Universidad Carlos III de Madrid, Spain; moscoso@math.uc3m.es} \and Alexei Novikov\footnote{Department of Mathematics,
Pennsylvania State University, USA; anovikov@math.psu.edu} \and George Papanicolaou\footnote{Department
of Mathematics, Stanford University, USA; papanicolaou@stanford.edu} \and Lenya Ryzhik\footnote{Department
of Mathematics, Stanford University, USA; ryzhik@math.stanford.edu}}

\begin{document}
\maketitle


\begin{abstract}
We present an ordinary differential equations approach to the
analysis of algorithms for constructing $l_1$ minimizing solutions
to underdetermined linear systems of full rank. It involves
a relaxed minimization problem whose minimum is independent of the
relaxation parameter. An advantage of using the ordinary differential equations
 is that  energy methods can be used to prove convergence. The 
connection to the discrete algorithms is provided by the Crandall-Liggett
theory of monotone nonlinear semigroups. We illustrate the effectiveness
of the discrete optimization algorithm in some sparse array imaging problems.
\end{abstract}

\section{ Introduction}

We consider the solution of large underdetermined linear systems of equations 
$Ax=y$ where $A\in \gR^{m \times n}$ is a given matrix, $y\in \gR^m$ is a
known vector of $m\ll n$ measurements, and $x\in \gR^n$ is the unknown
signal or image to be estimated. We assume that $A$ has full rank equal to $m$.
We want to find the solutions of this system with minimal $l_1$ norm $\|x\|_{l_1}$,
\begin{equation}
\min ||x||_{l_1},~\hbox{ subject to  } y = Ax \, . \label{l1}
\end{equation}
Our motivation is array imaging problems, which is an application 
discussed in this paper,  but such sparsity inducing
constrained minimization problems, where the $l_1$ norm of the
solution vector is used,  arise in many
other applications in signal and image processing \cite{Chan_book}.

%

A lot of research has been devoted to developing algorithms for solving
efficiently \eqref{l1} and its relaxed form
\begin{equation}
\min_x \lcu \tau ||x||_{l_1} + \frac{1}{2} \|y - Ax\|^2 \rcu\, .
\label{l1_noise}
\end{equation}
Here, and throughout the paper,
$\|q\|$ denotes the $l_2$-norm of a vector $q$. In
\eqref{l1_noise}, the exact constraint $y= A x$ is relaxed so as to take
into account possible measurement noise, and $\tau$ is a positive real
parameter that promotes sparse solutions when it is large enough.

The iterative shrinkage-thresholding algorithm (ISTA) is the usual 
gradient descent method applied to (\ref{l1_noise}).
It has been used in many
different applications with great success, such as
\cites{Daubechies05,Donoho95,Chambolle98,DeMol02,Figueiredo03,Yin08}, just to
mention a few. The ISTA algorithm generates a sequence of iterates 
$\{x_k \}$ of the form
\begin{equation}
x_{k+1} = \eta_{\tau h} (x_k - h\nabla f(x_k))\, . \label{ista}
\end{equation}
Here, $h$ is the step size,
\begin{equation}
\eta_a(x)=\begin{cases}
&x-a,~\hbox{ if } x>a,\\
&0,~\hbox{ if } -a< x<a,\\
&x+a,~\hbox{ if } x<-a\, \\
\end{cases}
\label{thresholding}
\end{equation}
is the shrinkage-thresholding operator, and $\nabla f(x_k)$
denotes the gradient of $f(x) = \frac{1}{2} ||y - Ax||^2$ at the
current iterate $x_k$. Thus, $\nabla f(x_k)=A^*(Ax_k - y)$, where
$A^*$ denotes the complex conjugate transpose of $A$. The algorithm
\eqref{ista} involves only simple matrix-vector multiplications
followed by a shrinkage-thresholding step.


For a fixed value of $\tau$ the
solution to \eqref{l1_noise} differs in general from the solution of
\eqref{l1}. In other words,  exact recovery from noiseless data
is not achieved by solving \eqref{l1_noise}, unless the regularization parameter
$\tau$ is sent to zero. However, it is well known that the
convergence of \eqref{ista} is slow for small values of
the parameter $\tau$. This issue is considered in detail in \cite{Bredies08}. 
Variants of \eqref{ista} have been
 proposed to speed up its convergence rate. In
\cite{Beck09}, for example, a fast version of ISTA is proposed (FISTA, described in more
detail below in Section~\ref{sec:numerics}) that has as easy
an implementation as \eqref{ista} but has a much better
convergence rate.

In this paper, we present an ordinary differential
equations (ODE) approach to an iterative
shrinkage-thresholding algorithm for solving
$\ell_1$-minimization problems independently of
the regularization parameter $\tau$.
We use a generalized Lagrange multiplier, or augmented Lagrangian,
approach \cite{bertsekas,FG_83,Hestenes,Powell,Rockafellar} to the relaxed problem \eqref{l1_noise} to impose exact
recovery of the solution to \eqref{l1}. The exact solution is
sought through an efficient algorithm obtained from a min-max
variational principle, which is a special case of the Arrow-Hurwitz-Uzawa
algorithm~\cite{AHU}. We prove that this algorithm
yields the exact solution for all values of the parameter $\tau$.
Our only assumption is that the matrix $A$ has full rank. 
The connection of the ODE method to the iterative shrinkage algorithm
is provided by the Crandall-Liggett theory \cite{CL}, which analyzes the convergence
of an implicit finite difference discretization of the ODE. The theory works for
infinite dimensional, monotone nonlinear problems as well.
The performance of the algorithm, with and without noise in the
data, is explored through several numerical simulations of array
imaging problems.

The min-max variational principle used here is also
behind the Bregman and linearized Bregman iterative
algorithms \cite{Goldsteis08,Osher05,Yin08,Yin}.
The fully implicit version of the algorithm is also analyzed in
detail in \cite{CP,EZC} using different techniques.
Many other methods
have been proposed in the literature to solve
\eqref{l1} and \eqref{l1_noise} with large data. We just mention here some of
them:
homotopy~\cites{Tibshirani94,Osborne00, Efron04}, interior-point methods~\cite{Wright05},
gradient
projection~\cite{Figueiredo07}, 
and proximal gradient in combination with
iterative shrinkage-thresholding~\cites{Nesterov83, Nesterov07,
Beck09}. A detailed discussion and analysis of 
monotone operator splitting methods can be found in~\cite{Fadili}.

Finding the constrained, minimal $l_1$ norm solution in (\ref{l1})
does not imply that this solution vector has minimal support, even
though the $l_1$ norm is sparsity promoting. Nevertheless in many
applications, in imaging in particular, this optimization method does
produce the minimal support, or minimal $l_0$ norm solution.
The theory of compressed sensing 
\cite{Donoho05,Candes05,Tsaig05,Candes06a, Donoho06,Candes06b,Candes06c} gives conditions
under which constrained $l_1$ and $l_0$ minimizations are equivalent.
We do not address this issue here.

The paper is organized as follows. In Section
\ref{sec:main} we motivate our approach, summarize our main
results, and describe the numerical algorithm. Theorems~\ref{thm-jan7} and
\ref{thm-jan10} are the main results of this paper. A key
ingredient in the proof of these theorems is Theorem~\ref{thm-ode}
proved in Section \ref{sec:proofs1}. The proof of the variational
principle of Theorem~\ref{thm-minimexact} is
presented in Section~\ref{sec:proofs2}.  This result is originally due
to~\cite{Rockafellar} but we present it here for the convenience of the reader.
In Section \ref{sec:numerics} we show
the performance of the algorithm with and without noise in the
data using some numerical experiments of array imaging. Finally,
Section \ref{sec:conclusions} contains conclusions.

{\bf Acknowledgment.} MM was supported by the Spanish Ministry of
Education and Science grant FIS2010-18473, AN was supported by
NSF grant DMS-0908011, LR was supported by AFOSR NSSEFF
fellowship, and NSF grant DMS-0908507, and GP by
AFOSR grant FA9550-11-1-0266. We thank Laurent Demanet
for an interesting discussion, Alberto Bressan for bringing
reference~\cite{CL} to our attention, and Jalal Fadili for taking
time to explain the literature on Arrow-Hurwicz-Uzawa algorithm, 
and pointing out the paper~\cite{CP} to us.

\section{Formulation and main results}
\label{sec:main}

We consider the constrained optimization problem (\ref{l1})
under the assumptions  that (\ref{l1}) has a
unique minimizer $\bar x$, and that $A$ has full rank:
the matrix $AA^*$ is invertible.

\subsection{The min-max variational principle}

In order to find the minimizer $\bar x$, we  recall  the variational
formulation of the $l_1$-minimization 
problem ~\cite{bertsekas,Hestenes,Powell,Rockafellar}.
Define the finction
\[
F(x,z)=\tau ||x||_{l_1} + \farc{1}{2}\|Ax-y\|^2+ \ltr z, y- Ax  \rtr,
\]
for $x\in\Rm^n$ and $z\in\Rm^m$, and set
\begin{equation}
\bar F=\max_z \min_x \lcu F(x,z) \rcu \, .
\label{gelma}
\end{equation}
\begin{prop}\label{lem1}
The problem (\ref{gelma}) has a solution, that is $-\infty<\bar F<+\infty$,
and the max-min is attained.
\end{prop}
{\bf Proof.}  The function $F(x,z)$ is convex in $x$, and
$\lim_{x\to\infty} F(x,z)=+\infty$, for any   fixed $z$. Thus, $F(x,z)$
attains its minimum for a fixed $z$.
Let us denote
\begin{equation}\label{nov1402}
l(x)=\tau ||x||_{l_1} + \farc{1}{2}\|Ax-y\|^2,
\end{equation}
and
\begin{equation}\label{lt}
h(z)= \min_{x} F(x,z)=\min_{x}[l(x)+\ltr z, y- Ax  \rtr] .
\end{equation}
As the function $l(x)$ is convex, and $l(x) \to +\infty$, as $|x| \to \infty$,
it follows 
that $h$ is  concave, as a minimum of affine functions,
and $h(z) \to -\infty$, as $|z| \to \infty$. Thus,
it attains its maximum
$\max_z h(z)$.~$\Box$

In order to motivate the functional~\eqref{gelma} we look at
another natural way to impose the constraint in~\eqref{l1} by using a Lagrange multiplier. If we consider a functional
\begin{equation}\label{degenerate}
\tau ||x||_{l_1}+\langle z, y-A x \rangle,
\end{equation}
then (at least, formally) its Euler-Lagrange equations  for the extremum
give us the {\it sub-differential} optimality condition
\begin{equation}\label{main}
\left[ A^*z \right]_i= \begin{cases}
&\tau, \hbox{ if } \bar{x}_i >0,\\
&-\tau, \hbox{ if } \bar{x}_i <0,\\
\end{cases}
\hbox{ and } \left| \left[ A^*z \right]_i\right| \leq \tau.
\end{equation}
It is, however, difficult to work with~\eqref{degenerate},
because if some of the entries of $A^*z$ are larger than $\tau$ in absolute value, then~\eqref{degenerate}
 is not bounded from below as a function of $x$. Further, even if $z$ is chosen
according to the sub-differential condition~\eqref{main}, then the minimum may not be unique, even if $A$ is invertible. Indeed, consider a simple example: minimize $|x|$ subject to $x=1$.
Suppose $\tau=1$, then~\eqref{degenerate} is $|x| +z(1-x)$. Then $z=1$ satisfies the sub-differential condition,
and~\eqref{degenerate} becomes
\[
|x| +(1-x)=\begin{cases}
&1, \hbox{ if } x >0,\\
&1- 2 x, \hbox{ if } x <0,\\
\end{cases}
\]
which has no minimum.
The addition of a quadratic term to (\ref{degenerate}) regularizes this degeneracy. 
Since the function $l(x)$ in~\eqref{nov1402}
is convex, \eqref{lt} may be interpreted (up to a sign) as a generalized Legendre transform of $l(x)$. 

The first observation is 
that if (\ref{l1})  has a unique minimum $\bar x$ then the variational principle
(\ref{gelma}) finds $\bar x$ exactly.
\begin{thm} \label{thm-minimexact}
Assume that  (\ref{l1})  has a unique minimum $\bar x$.
Then we have
\begin{equation}\label{minmax}
\tau ||\bar{x}||_{l_1}= \max_{z} \min_{x} F(x,z)
\end{equation}
Moreover, we have $\tau ||\bar{x}||_{l_1}=F(\bar{x},z)$ for any $z$,
and if  $\min_{x} F(x,z) = \tau ||\bar{x} ||_{l_1}$ for some fixed $z$, then $\argmin_{x} F(x,z)=\bar{x}$.
\end{thm}
This result can be found in~\cite{Rockafellar} in a much greater
generality. We present
its proof below in the particular case we are interested in,
for convenience of the reader.

It is remarkable that (\ref{minmax}) holds for any
value of $\tau>0$ -- this gives us a freedom to choose $\tau$ large or small, depending
on a particular application.
We also have the following well known result \cite{Rockafellar}, which 
follows from the proof of Theorem~\ref{thm-jan7} below.
\begin{thm}\label{thm-sub}
Assume that (\ref{l1}) has a unique minimizer $\bar x$.
Then, there exists a vector $z$ such that
$[A^*z]_i=\hbox{sgn}(\bar x) $ if $\bar x\neq 0$, and
$|[A^*z]_i|\le 1$ if $\bar x_i=0$.
\end{thm}
We say that $z$ satisfies the sub-differential condition if there exists a minimizer of (\ref{l1}) such that
\begin{equation}\label{nov1006}
\hbox{$[A^*z]_i=\tau\hbox{sgn}(\bar x) $ if $\bar x\neq 0$, and $|[A^*z]_i|\le \tau$ if $\bar x_i=0$.}
\end{equation}
We note that (\ref{nov1006}) is weaker than the sub-differential condition of~\cite{Candes06a} --
there it is required that  $|[A^*z]_i|< \tau$ if $\bar x_i=0$,
while we do not require the strict inequality.
It follows from the proof of Theorem~\ref{thm-sub} that the exact extremum of $F(x,z)$ is achieved for { any} $z$ that satisfies the sub-differential condition~\eqref{nov1006}.
Such $z$ is not unique but, of course, our interest is not in finding $z$ but in finding  the minimizer of (\ref{l1}).

\subsection{The ordinary differential equations method}

In order to find $\bar x$, ideally,
we would like to take the ODE point of view and
generate a trajectory $(x(t), z(t))$ of the following system
\begin{eqnarray}
\frac{dx}{dt} = - \nabla_x F(x,z),~~
\frac{dz}{dt} =  \nabla_z F(x,z),  \label{subgradsch2bis}
\end{eqnarray}
with the hope that $x(t)\to \bar x$ as $t\to +\infty$. There is an obvious degeneracy
in the problem, namely, $F(\bar x,z)=\tau\|\bar x\|_{l_1}$ for all $z\in\Rm^m$.
Hence, we can only hope to recover $\bar x$ as there
is no "optimal" $z$.

The obvious technical difficulty is that
the function $F(x,z)$ is not differentiable in $x$ at the points where $x_j=0$ for some $j=1,\dots,n$.
Following~\cite{CL}, we interpret solutions of (\ref{subgradsch2bis}) as follows. Given $x\in\Rm^n$, let
the sub-differential $\partial\|x\|_{l_1}$ be a subset of $\Rm^n$:
\[
\partial\|x\|_{l_1}=\sgn(x_1)\times\dots\times \sgn(x_n).
\]
Here $\sgn(s)$, for $s\in\Rm$, is understood a subset
of $\Rm$: $\sgn(s)=\{1\}$ if $s>0$, $\sgn(s)=\{-1\}$ if $s<0$ and $\sgn(s)=[-1,1]$ if $s=0$. Then, instead of treating
the system of ODEs (\ref{subgradsch2bis}) with a discontinuous right side, we consider
\begin{eqnarray}\label{jan602}
&&\farc{dx}{dt}-A^*(z-Ax+y)\in -\tau\partial\|x\|_{l_1},\\
&&\farc{dz}{dt}=y-Ax,\nonumber
\end{eqnarray}
supplemented by the initial data $x(0)=x_0$, $z(0)=0$. We say that
$(x(t),z(t))$ is a strong solution to (\ref{jan602}) on a time
interval $0\le t\le T$ if $x(t)$ and $z(t)$ are continuous,
differentiable for almost all $t\in[0,T]$, $x(0)=x_0$, $z(0)=0$,
and (\ref{jan602}) holds for almost all $t\in[0,T]$. 

An important observation is that (\ref{jan602}) 
is contractive, or, accretive in the sense of Crandall and Liggett~\cite{CL}.
That is, the following property holds: given any pair $(x_1,z_1)$, $(x_2,z_2)$ and any $\xi_1\in\partial\|x_1\|_{l_1}$,
$\xi_2\in\partial\|x_2\|_{l_1}$, we have:
\begin{eqnarray}
&&(A^*(z_1-Ax_1)-\tau \xi_1-A^*(z_2-Ax_2)+\tau\xi_2)\cdot(x_1-x_2)-(Ax_1-Ax_2)\cdot(z_1-z_2)\nonumber\\
&&=-\tau(\xi_1-\xi_2)\cdot(x_1-x_2)- \|A(x_1 -x_2)\|^2 \le 0.\label{jan1804}
\end{eqnarray}
The last inequality above follows from the component-wise monotonicity of the sub-differential $\partial\|x\|_{l_1}$. It follows from (\ref{jan1804}) and
Theorems I and II of~\cite{CL} that (\ref{jan602}) has a unique strong solution.
Our first result shows that this solution converges as $t\to+\infty$ to 
$\bar x$, the minimizer of (\ref{l1}).
\begin{thm}\label{thm-jan10}
Let (\ref{l1}) have a unique minimizer $\bar x$. Then, for
any $\delta>0$ there exists $T=T(\delta)$ such that the solution
of~\eqref{jan602}~ satisfies
\begin{equation}\label{jan1012}
\|x(t) -\bar{x}\|< \delta,\hbox{ for all $t>T$.}
\end{equation}
The time $T(\delta)$ depends only on $\delta$, the initial data
$x_0$, and $\|AA^*\|$ but not on the dimension $n$.
\end{thm}

\subsection{The discrete algorithm}

 We consider the following numerical algorithm to solve (\ref{jan602}):
\begin{eqnarray}
&&\frac{x_{k+1} - x_k}{\Delta t} = -\tau\xi_{k+1}+A^*(z_{k+1}+y-Ax_{k+1})
\label{implicit1} \, ,\\
&&\frac{z_{k+1} - z_k}{\Delta t}=y-Ax_{k+1},   \nonumber \,
\end{eqnarray}
with the initial data $x_0=x$, $z_0=0$. Here, $\xi_{k+1}$ is a
vector in the set $\partial\|x_{k+1}\|_{l_1}$.

A simple way to understand how (\ref{implicit1}) works is to consider
the toy problem
\begin{equation}\label{jan704}
\dot r=-\sgn{r}.
\end{equation}
An explicit discretization
\[
\farc{r_{k+1}-r_k}{\Delta t}=-\xi_{k},
\]
with $\xi_{k}\in\sgn(r_{k})$, will start oscillating around $r=0$
as soon as $r_k\in[-\Delta t,\Delta t]$, and will never converge
to $x=0$ for $\Delta t>0$. On the other hand, the implicit
discretization
\begin{equation}\label{jan706}
\farc{r_{k+1}-r_k}{\Delta t}=-\xi_{k+1},
\end{equation}
with $\xi_{k+1}\in\sgn(r_{k+1})$ behaves differently.
If $r_k\in[-\Delta t,\Delta t]$, the implicit nature of this scheme shows that
it is impossible to have $\xi_{k+1}=\pm 1$, which 
forces $\xi_{k+1}=r_k/\Delta t$ and
$r_{k+1}=0$. The implicit scheme 
is actually equivalent to soft thresholding:
\begin{equation}\label{jan708}
r_{k+1}=\eta_{\Delta t}(r_k).
\end{equation}
The function $\eta_s$ here is defined by (\ref{thresholding}).
This simple example already shows both 
the importance of using an implicit discretization, and that the implicit
scheme has a simple explicit realization (\ref{jan708}).

Theorems I and II of~\cite{CL} not only provide existence of a strong
solution to (\ref{jan602}) but also show that it can be found by the 
implicit scheme (\ref{implicit1}).
\begin{prop}\label{prop-jan6}
Solution of (\ref{implicit1}) converges as $\Delta t \to 0$, uniformly on finite time intervals, to the unique strong solution of (\ref{jan602}).
\end{prop}
Theorem~\ref{thm-jan10} and Proposition~\ref{prop-jan6}
together imply immediately the following theorem.
\begin{thm}\label{thm-jan7}
Let the sequence $x_n$, $z_n$ solve (\ref{implicit1}) with the
initial data $x_0=x$, $z_0=0$. Given any $\delta>0$ there exists
$h>0$ and $T>0$, so that for all $0<\Delta t<h$ and all
$k>[T/\Delta t]$ we have $|x_k-\bar x|<\delta$. The time $T$
depends on $\delta$, the initial data $x\in\Rm^n$, and the norm
$\|AA^*\|$.
\end{thm}

If one examines the proof of Theorems I and II in~\cite{CL}, it is clear that the
only term that should be discretized implicitly is $\hbox{sgn} x$ -- the other terms
can be discretized explicitly, keeping the statement of Proposition~\ref{prop-jan6}
intact. Hence,
the result  of Theorem~\ref{thm-jan7} applies equally well
to an Euler quazi-explicit  modification
of (\ref{implicit1}) that is
easier to implement numerically:
\begin{eqnarray}
&&x_{k+1} = x_k -   \xi_{k+1} +\Delta t A^*(z_k+y-Ax_k)
\label{fd1bis} \, ,\nonumber\\
&&z_{k+1} =  z_k + \Delta t  (y - A x_k) \label{fd2bis},
\end{eqnarray}
where $\xi_{k+1} \in \tau \Delta t \,\partial || x_{k+1} ||_{l_1}$
is a vector in the subdifferential of $\tau \Delta t\, || x_{k+1}
||_{l_1}$. We call this scheme the generalized Lagrangian
multiplier algorithm (GeLMA). As in the toy problem
(\ref{jan704})-(\ref{jan708}), it is equivalent to soft
thresholding:
\begin{eqnarray}
&&x_{k+1}=\eta_{\tau \Delta t} \lc x_k +\Delta t A^*(z_k+y-Ax_k)\rc
\nonumber \, ,\\
&&z_{k+1} =  z_k + \Delta t  (y - A x_k) \label{fd2bis2}.
\end{eqnarray}
This scheme converges if $\Delta t<1/||A||$ -- that condition simply
comes from the usual constraint for an explicit scheme for a linear system.
GeLMA algorithm is extremely easy to implement numerically.

We also note that one can mimic  the ODE proof of Theorem~\ref{thm-jan10} directly
on the numerical scheme, eliminating, in particular, the dependence of $h$ on $\delta$.
Our objective, however, in part, is to explain the effectiveness of
shrinking-thresholding algorithms in the language of differential equations,
potentially opening the way for the application of other continuous techniques in such
problems. Therefore, we have chosen to concentrate on the ODE proof.




%


\subsubsection{The regularized ordinary differential equations}

Since the system (\ref{jan602}) has a "bad" right side, working with it 
directly is technically inconvenient. Hence, in order to prove
Theorem~\ref{thm-jan10}, from which Theorem~\ref{thm-jan7}
follows, we consider a regularized system,
%
introducing a single-valued approximation of $\hbox{sgn} x$:
\[
G_\eps (s)= \begin{cases}
&1,~\hbox{ if } s>\eps,\\
&s/\eps,~\hbox{ if }  |s|<\eps,\\
&-1,~\hbox{ if } s<-\eps.\\
\end{cases}
\]
Here $\eps>0$ is a small regularization parameter that will be sent to zero
at the end. With a slight abuse of notation, here,
and in other instances when this should cause no confusion, we will also denote by $G_\eps(x)$ a vector valued  
function  with components $G_\eps(x)=(G_\eps(x_1),G_\eps(x_2),\dots,G_\eps(x_n))$.
The regularized version of (\ref{jan602}) is
\begin{eqnarray}
\frac{dx_\eps}{dt} = -\tau G_{\eps}\left(x_\eps \right) +A^* \left( z  +y-Ax_\eps \right) \label{s11bis},~~~~
\frac{dz_\eps}{dt} =  y - Ax_\eps.
\end{eqnarray}
It has the same form (\ref{subgradsch2bis}), with $F(x,z)$ replaced by a
differentiable approximation
\begin{equation}\label{nov804}
F_\eps(x,z)=\tau \sum_{j=1}^nr_\eps(x_j) + f(x) + \langle z, y- Ax  \rangle .
\end{equation}
Here,
\[
r_\eps(s)= \begin{cases}
&|s|,~\hbox{ if } s>\eps,\\
&s^2/(2\eps)+\eps/2,~\hbox{ if }  |s|<\eps,\\
&|s|,~\hbox{ if } s<-\eps,\\
\end{cases}
\]
is an approximation of $|s|$ known as the Huber function. 
We will denote below
\begin{equation}\label{nov826}
\|x\|_{l_\eps^1}=\sum_{j=1}^n r_\eps(x_j),
\end{equation}
though, of course, this is not a norm as it does not vanish at $x=0$.
\begin{thm}\label{thm-ode}
Let (\ref{l1}) have a unique minimizer $\bar x$. Then, for
any $\delta>0$ there exists $\eps_0=\eps_0(\delta,n)$ and
$T=T(\delta)$ such that for any $\eps$, $0<\eps<\eps_0$ the
solution of~\eqref{s11bis}~ satisfies
\begin{equation}\label{nov914}
\|x_\eps(t) -\bar{x}\|< \delta,\hbox{ for all $t>T$.}
\end{equation}
The time $T(\delta)$ depends only on $\delta$, the initial data $x_0$, and $\|AA^*\|$ but not on the dimension $n$.
\end{thm}
When the minimizer of (\ref{l1}) is not unique, the proof of Theorem~\ref{thm-ode}
can be easily adapted   to show that for any $\delta>0$ there exists $\eps_0(\delta)$ such that for any $\eps\in(0,\eps_0)$
and any limit point (as $t\to+\infty$)  $\bar x_\eps$ of the trajectory $x_\eps(t)$, we have $\|\bar x_\eps-\bar x\|<\delta$ for
some minimizer $\bar x$ of (\ref{l1}).

Theorem~\ref{thm-ode} is the key ingredient in the proof of Theorem~\ref{thm-jan7}:
together with a priori bounds on $x_\eps(t)$ obtained in the course of its proof, they
show that solution $x(t)$ of (\ref{jan602}) is the limit of $x_\eps(t)$ as
$\eps\to 0$, and thus it obeys the same bounds as $x_\eps(t)$, finishing the proof.

\section{Application to array imaging}
\label{sec:numerics}

In this section we illustrate the performance of our algorithm for
array imaging of localized scatterers. The problem is to determine
the location and reflectivities of small scatterers by sending a
narrow band (single frequency) probing signal of wavelength
$\lambda$ from an active array and recording the backscattered
field on this array \cite{Borcea02}. In this paper we consider
only single illumination by the central element of the array.

\subsection{Array imaging in homogeneous media}

The array has $N$ transducers located at positions $\vect x_p$ ($p
= 1, \dots ,N$) separated from each other by a given distance. In
each numerical experiment there are $M$ point-like scatterers of
unknown reflectivities $\rho_j>0$ located at unknown positions
$\vect y_{n_j}$ ($j=1,\dots,M$). The scatterers are assumed to be
within a bounded region at a distance $L$ from the array, called
the Image Window (IW). We discretize this IW with a uniform mesh
of $K$ points $\vect y_j$ ($j=1,\dots,K$), and assume that each
scatterer is located at one of these $K$ grid points, so $\{\vect
y_{n_1},\ldots,\vect y_{n_M}\}\subset\{\vect y_1,\ldots,\vect
y_K\}$.

Furthermore, we assume that the medium between the array and the
scatterers is homogeneous so wave propagation between any two
points $\vect x$ and $\vect y$ is modeled by the free space Green
function
\begin{equation}\label{greenfunc}
\wG_0(\vect y,\vect x,\omega)=\frac{\exp(-i \kappa|\vect x-\vect
y|)}{4\pi|\vect x-\vect y|}\, ,
\end{equation}
where $\kappa=\omega/c=2\pi/\lambda$, and $c$ is the reference
wave speed in the medium. We also assume that the scatterers are
well separated or are weak, so multiple scattering among them is
negligible (this is the Born approximation). Under these conditions, the
backscattered field measured at $\vect x_r$ due to a pulse sent
from $\vect x_s$, and reflected by the $M$ scatterers in the IW,
is given by
\begin{equation}\label{discretedata}
b_r(\omega)=\sum_{j=1}^M\rho_j\wG_0(\vect x_r,\vect
y_{n_j},\omega)\widehat G_0(\vect y_{n_j},\vect x_s,\omega)\, .
\end{equation}

Next, we write the linear system that relates the reflectivity
$\rho_{0j}$ at each grid point $\vect y_j$ of the IW
($j=1,\dots,K$) and the data $b_r(\omega)$ measured at the array
($r=1,\dots,N$). To this end, we introduce the {\em reflectivity
vector} $\bfrho_0 =
(\rho_{01},\rho_{02},\dots,\rho_{0K})^T\in\mR^K$ and the data
vector $\vect b(\omega) = (b_{1},b_{2},\dots,b_{N})^T\in\mR^N$,
where the superscript $T$ means transpose. Thus, the image is a
gridded array of $K$ pixels, and the data is stacked into a vector
of $N\ll K$ components. Furthermore, there are only a few
scatterers in the IW so the vector $\bfrho_0$ is sparse.

Let us consider the vector
\[
\vect \wg_0(\vect y_j,\omega)=( \wG_0(\vect x_{1},\vect
y_j,\omega), \wG_0(\vect x_{2},\vect y_j,\omega), \dots,
\wG_0(\vect x_{N},\vect y_j,\omega) )^T\, ,
\]
that represents the signal at the array due to a point source at
$\vect y_j$ in the IW. Due to the spatial reciprocity $\wG_0(\vect
x_{i},\vect y,\omega)=\wG_0(\vect y,\vect x_{i},\omega)$, it can
also be interpreted as the illumination vector of the array at
position $\vect y_j$. With this notation, we can write the linear
system
\begin{equation}
A_\omega \bfrho_0= \vect b(\omega) \, , \label{lsystem}
\end{equation}
where $A_\omega$ is an $N \times K$ matrix whose $j^\mathrm{th}$
column is given by
$\widehat G_0(\vect y_{j},\vect x_s,\omega)\, \vect\wg_0(\vect
y_j,\omega)$. Since $N\ll K$, \eqref{lsystem} is an
underdetermined linear system, and hence there can be many
configurations of scatterers that match the data vector $\vect
b(\omega)$. Array imaging  is to solve \eqref{lsystem} for
$\bfrho_0$.

A related problem to \eqref{lsystem} has been studied in
\cite{Chai11} in array imaging of localized scatterers from
intensity-only measurements. Intensity measurements are
interpreted as linear measurements of a rank one matrix associated
with the unknown reflectivities. Since the rank minimization
problem is NP-hard, it is replaced by the minimization of the
nuclear norm of the decision matrix. This makes the problem convex
and solvable in polynomial time. It is shown that exact recovery
can be achieved by minimizing this problem.

\subsection{Numerical Simulations}

We consider here numerical experiments in 2D. Our linear array
consists of $100$ transducers that are one wavelength apart.
Hence, the aperture of the array is $a=100$. In each numerical
experiment there are a few point-like scatterers of different
reflectivity at a distance $120$ from the array. The IW is
discretized with $41 \times 41$ grid points. Hence, we have $1681$
unknowns and $100$ measurements. All the spatial units are
expressed in units of the wavelength $\lambda$ of the
illuminating source.

\begin{figure}[htbp]
\begin{center}
\begin{tabular}{ccc}
\includegraphics[scale=0.25]{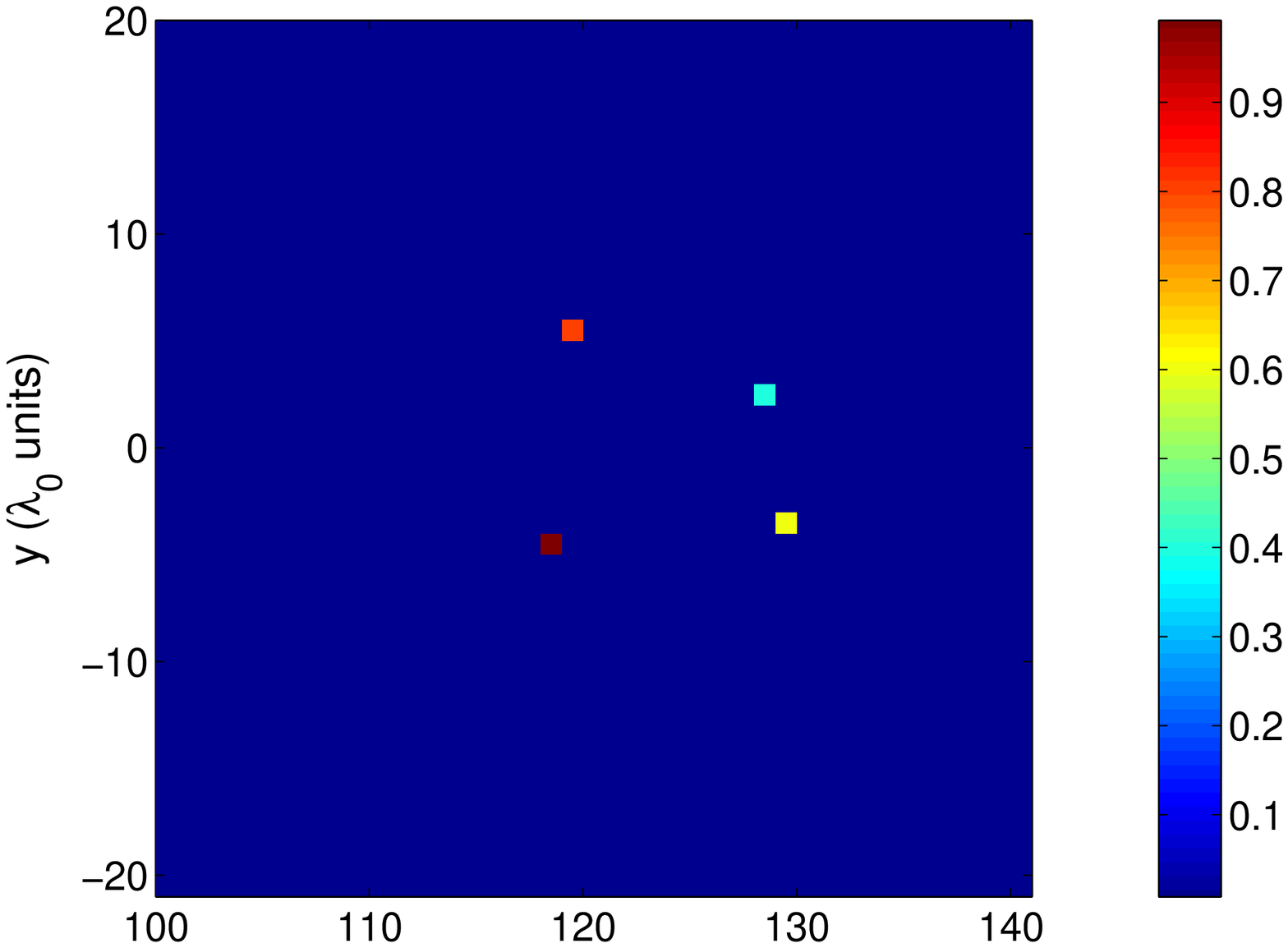}&
\includegraphics[scale=0.25]{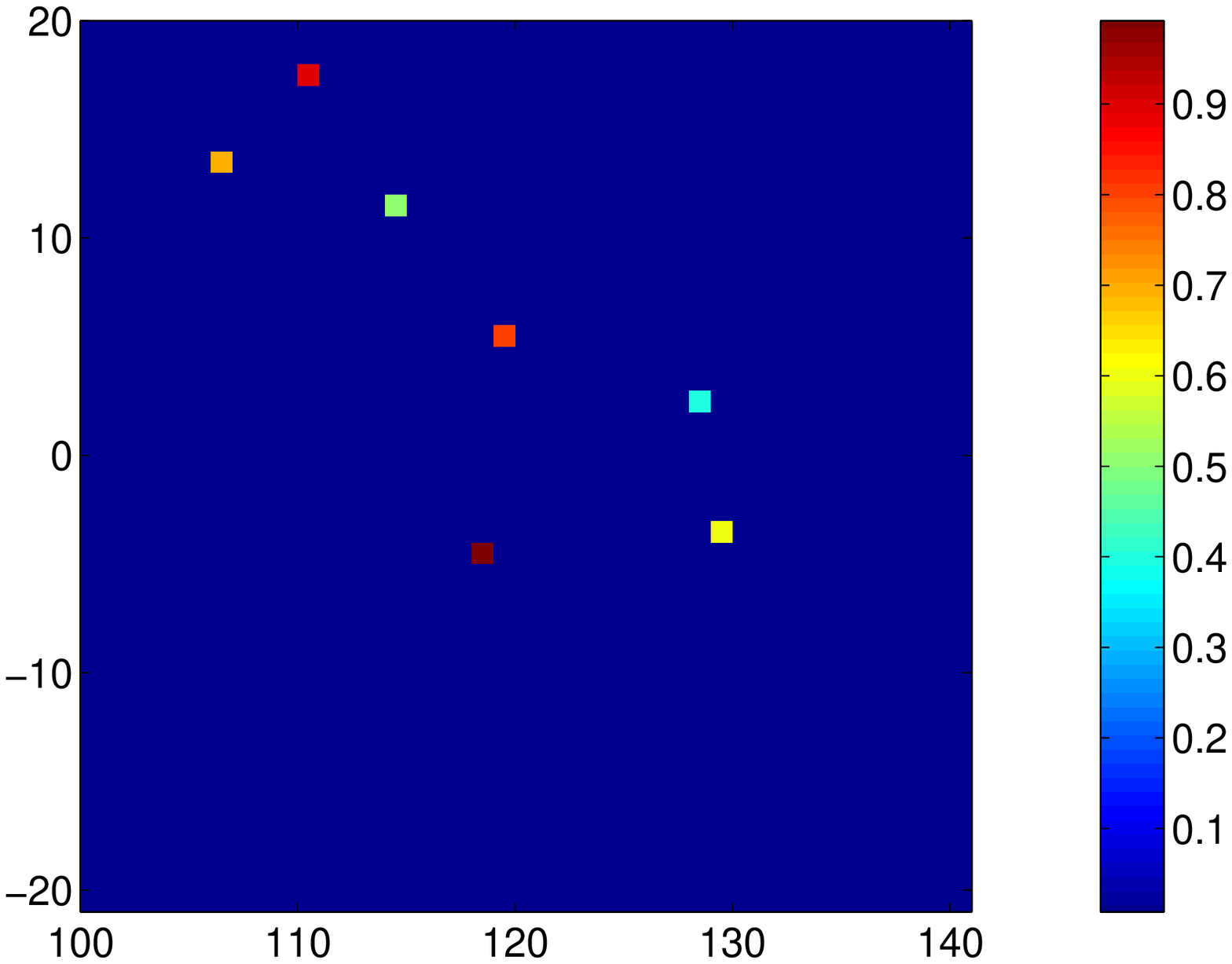}&
\includegraphics[scale=0.25]{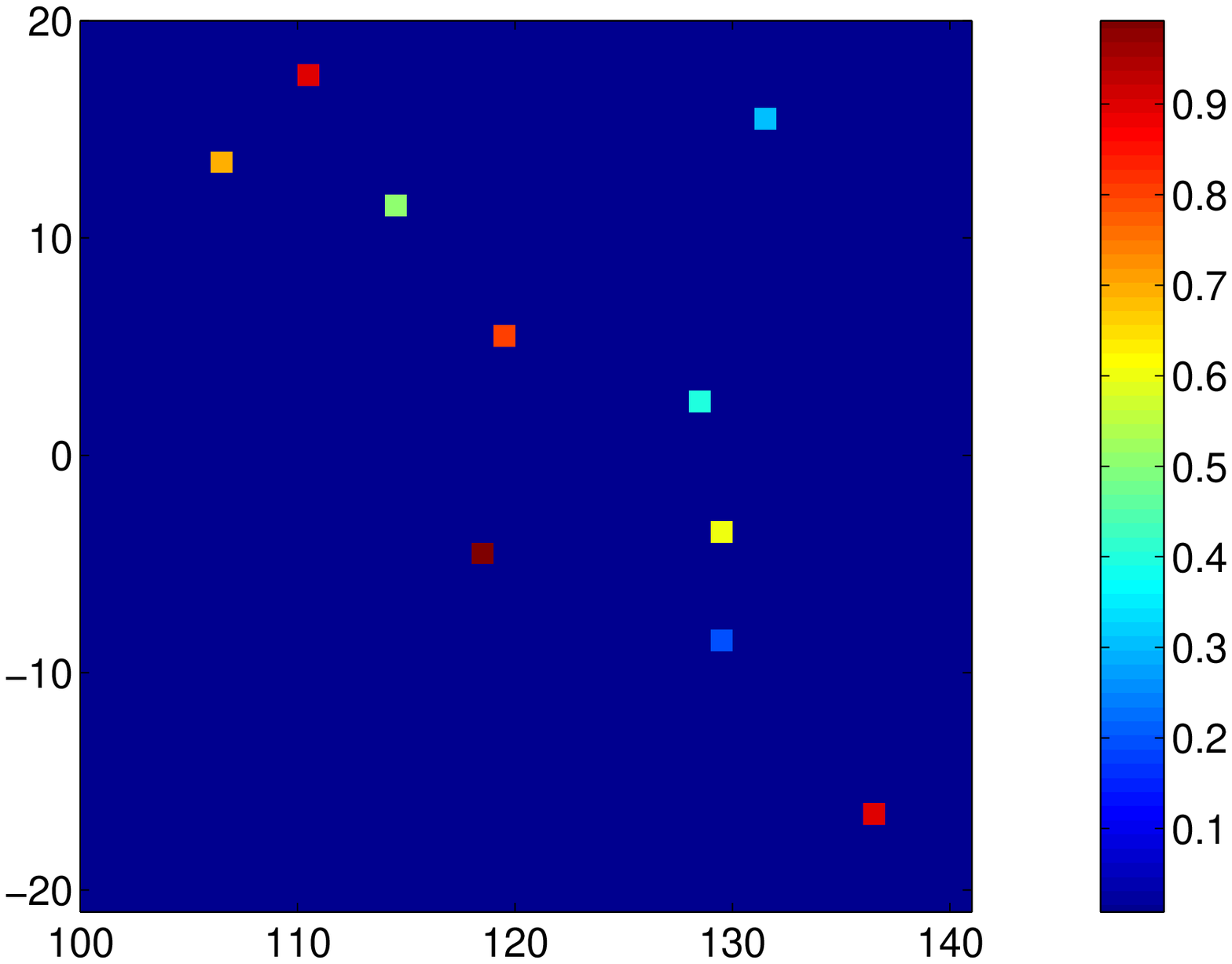}\\
\includegraphics[scale=0.25]{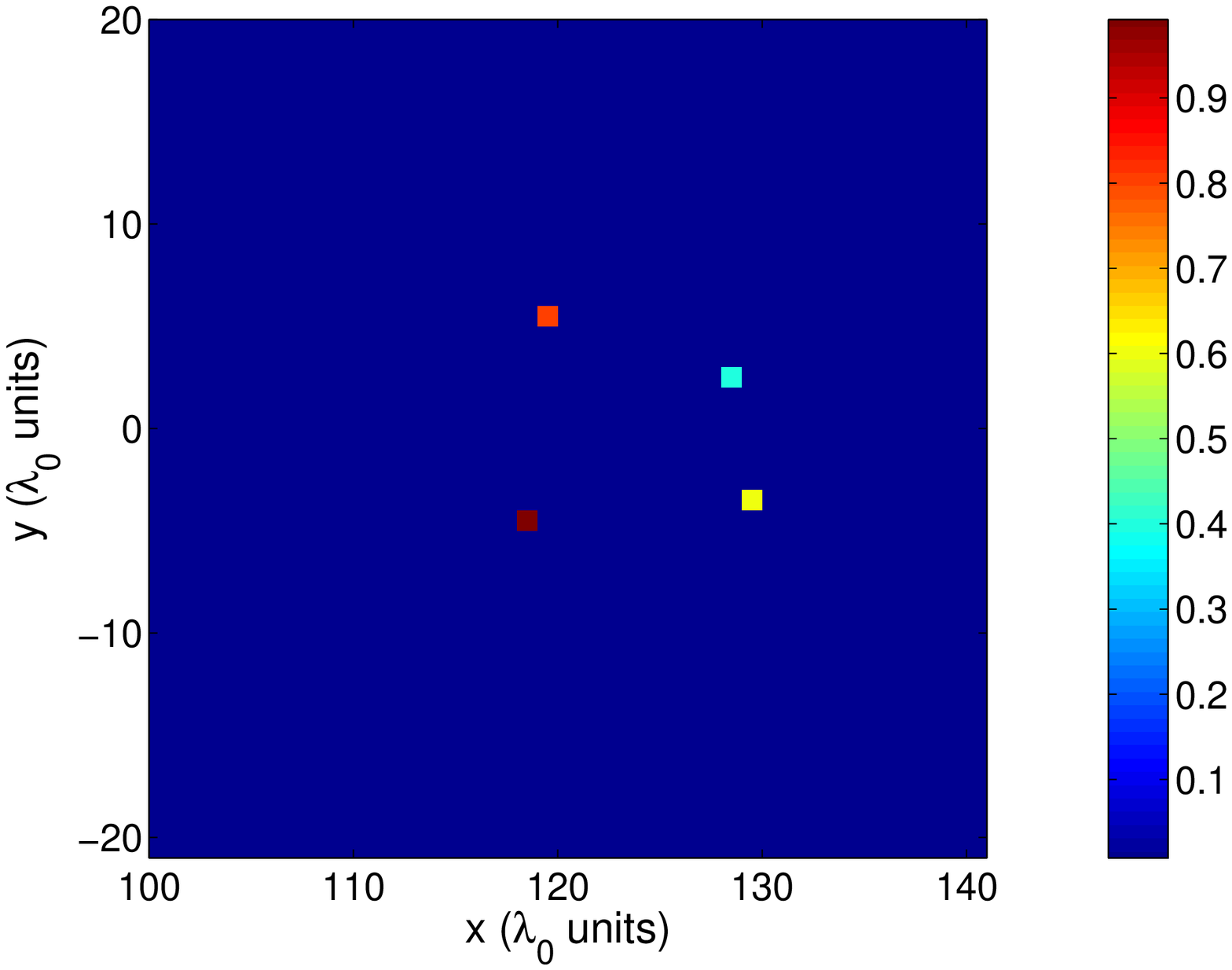}&
\includegraphics[scale=0.25]{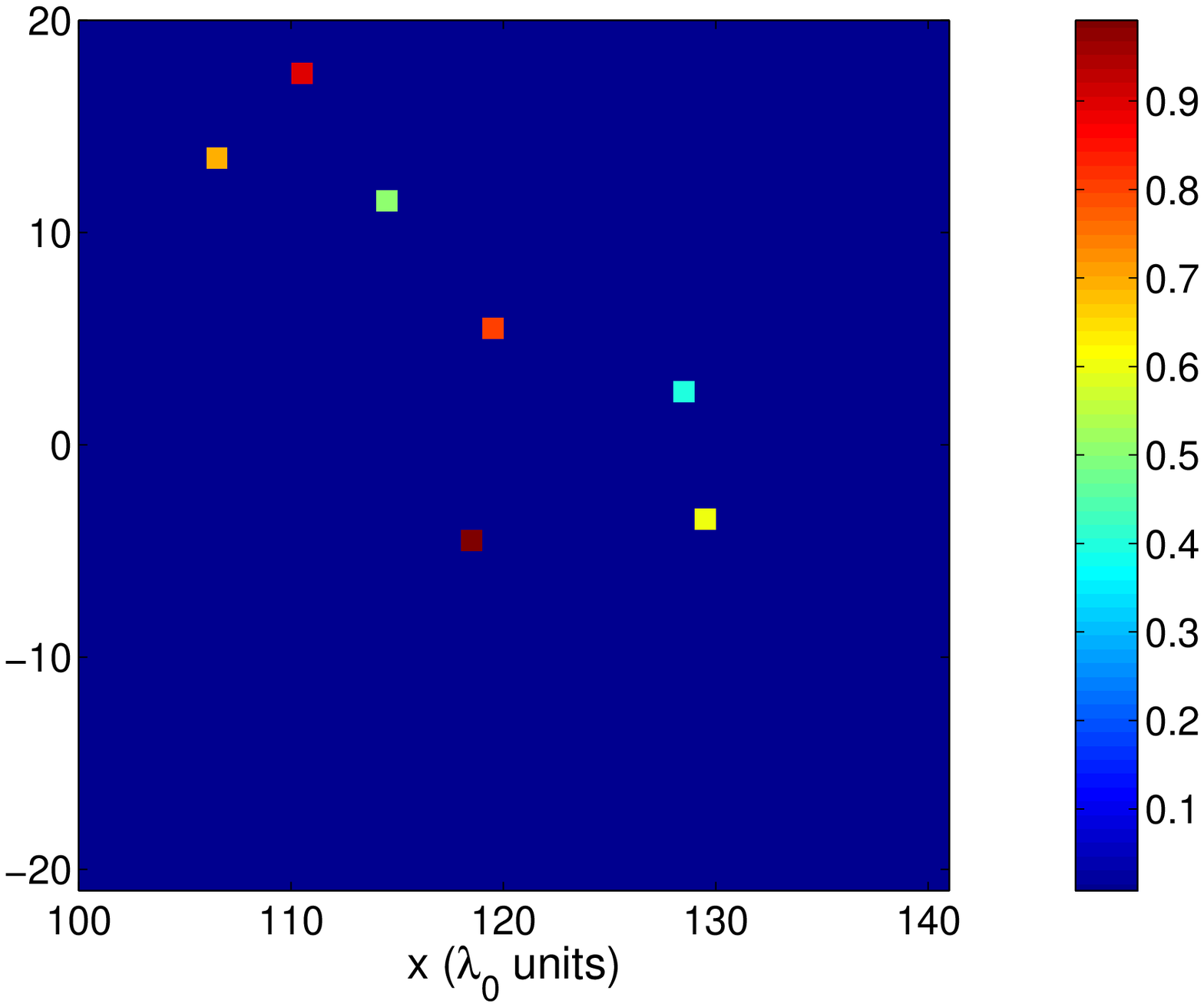}&
\includegraphics[scale=0.25]{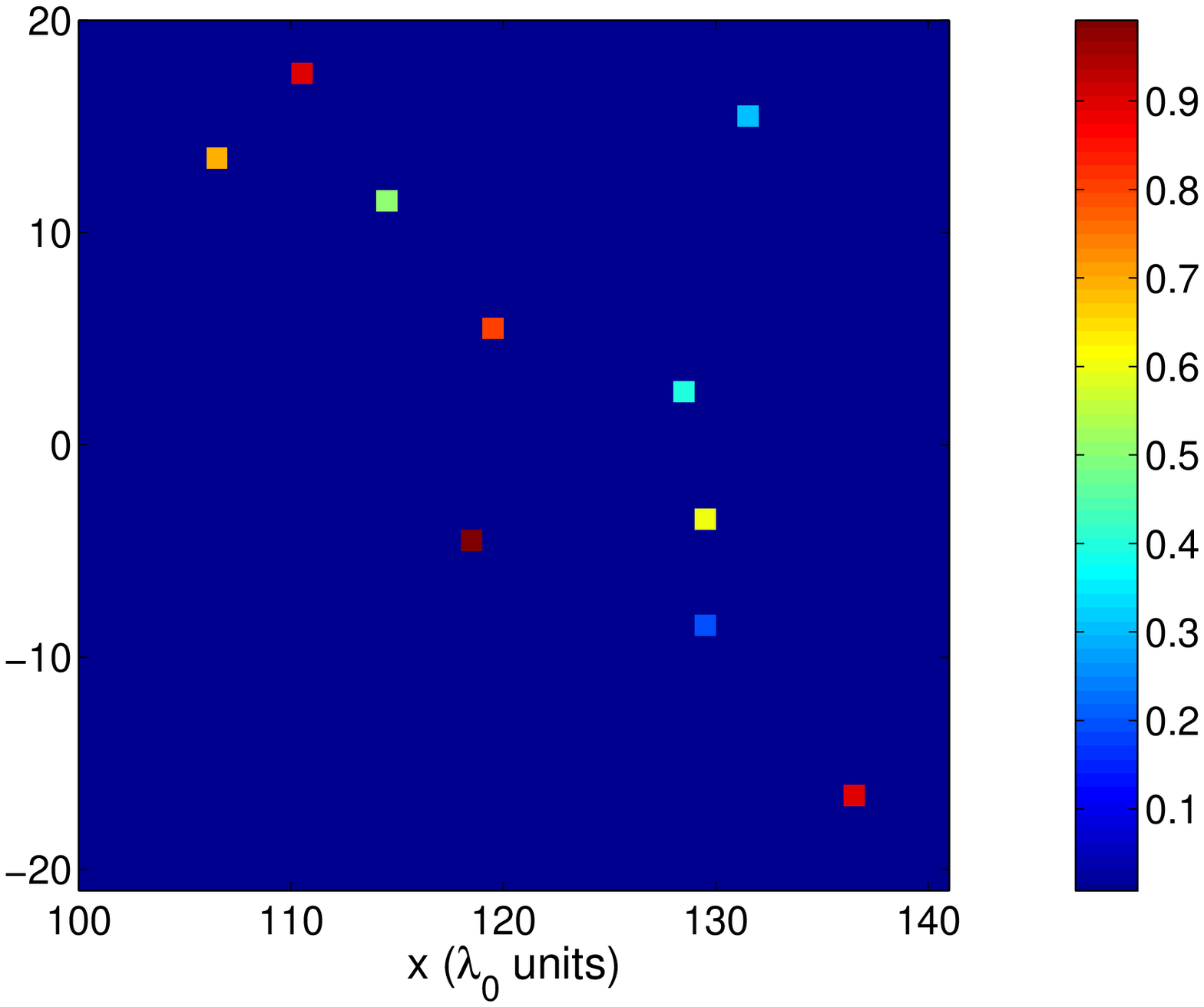}
\end{tabular}
\caption{Top row: original configurations of the scatterers within
the $41\times 41$ IW. Bottom row: recovered images obtained by the
$\ell_1$ minimization GeLMA algorithm \eqref{fd2bis2} with no
noise in the data.$\tau=20\|A_\omega^T{\bf
b}(\omega)\|_{l_\infty}$ and $N_{iter}=300$.}
\label{fig:imagesnonoise}
\end{center}
\end{figure}

\begin{figure}[htbp]
\begin{center}
\begin{tabular}{ccc}
\includegraphics[scale=0.25]{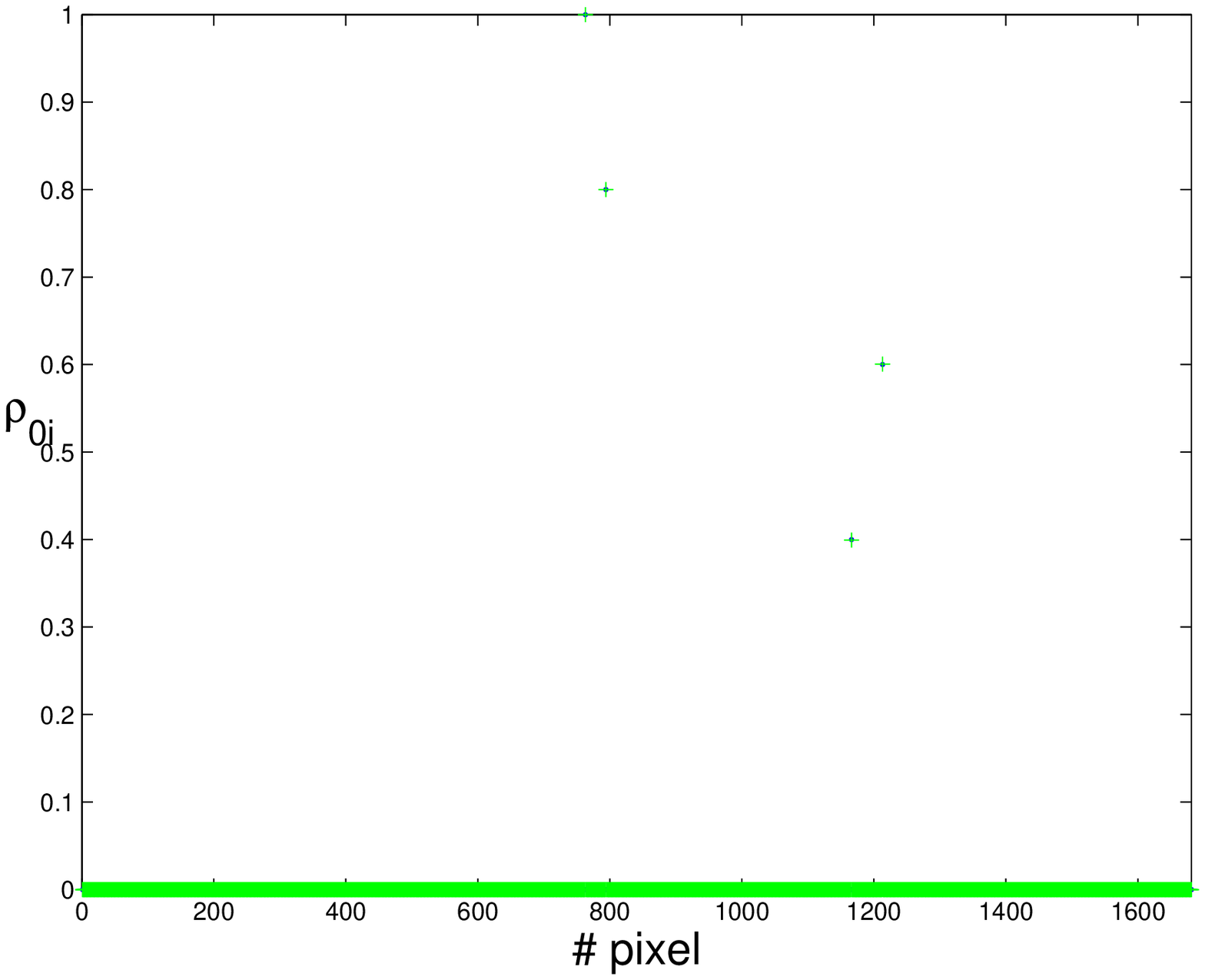}&
\includegraphics[scale=0.25]{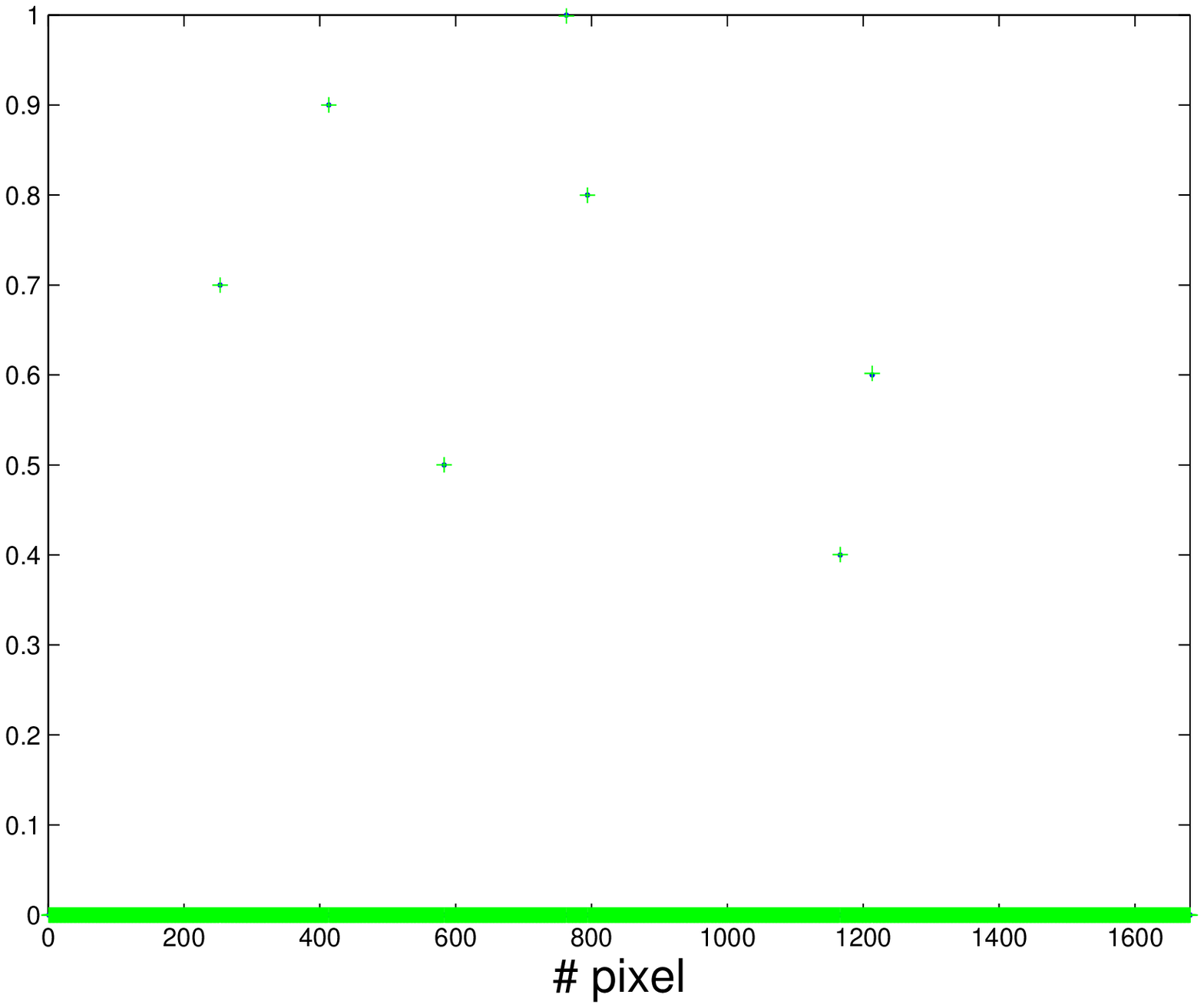}&
\includegraphics[scale=0.25]{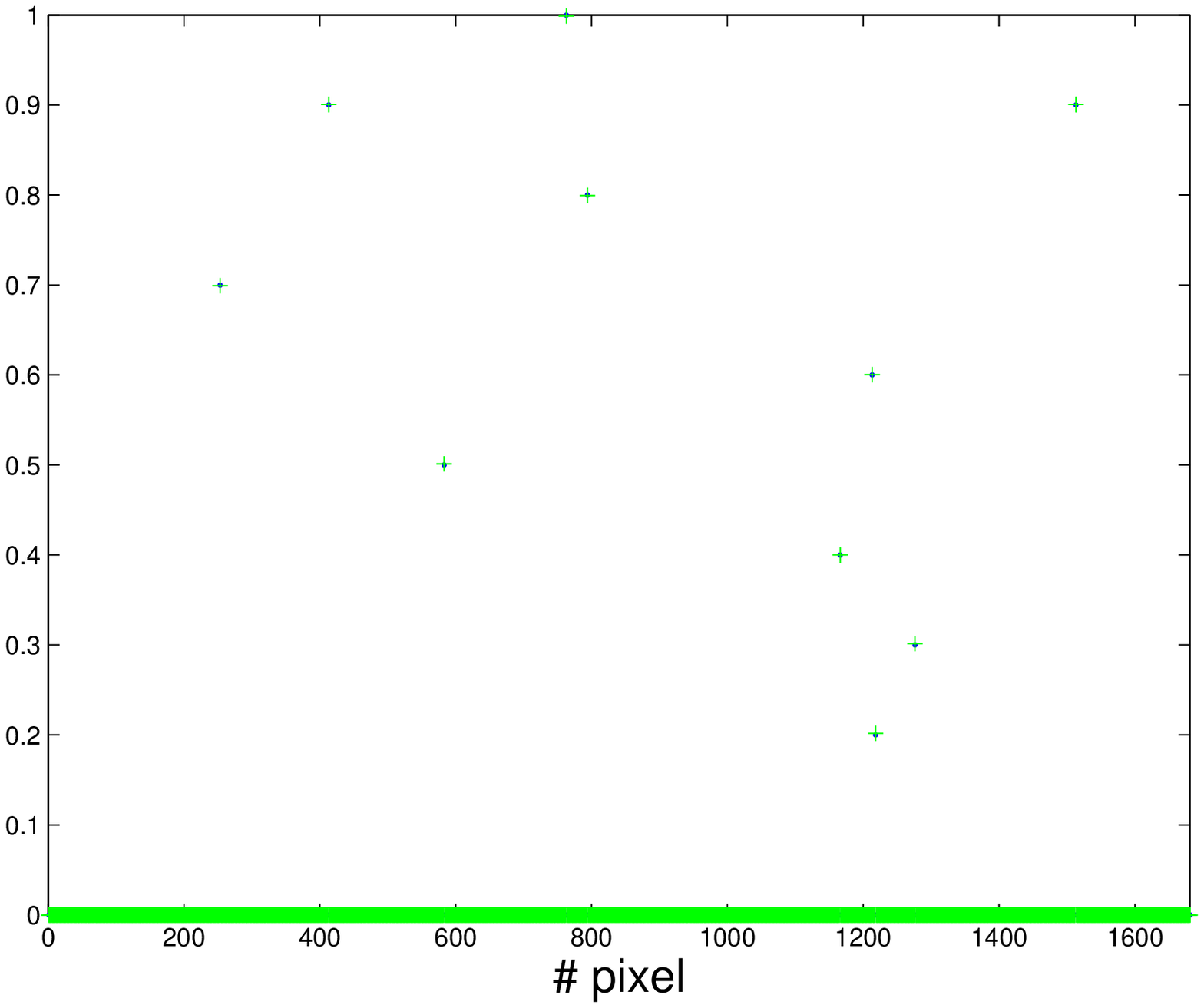}
\end{tabular}
\caption{Comparison between the exact solutions (blue circles) and
the solutions obtained with the GeLMA algorithm (green crosses)
with no noise in the data.} \label{fig:pixelnonoise}
\end{center}
\end{figure}

Fig. \ref{fig:imagesnonoise} shows results from various
scatterer's configurations with no noise in the data. In the top
row we display the original scatterer's configurations and in the
bottom row the corresponding images obtained by the $\ell_1$
minimization GeLMA algorithm \eqref{fd2bis2}. These results show
that this algorithm recovers the positions and reflectivities of
the scatterers exactly  when there is no noise in the data. To
examine this issue more clearly we plot in Fig.
\ref{fig:pixelnonoise} the vector solutions $\bfrho$ (green
crosses) and the exact vectors $\bfrho_0$ (blue circles) for these
three scatterer's configurations. There is not apparent difference
between the exact and recovered solutions. Both, localization
(support recovery) and strength estimation (reflectivities) are
solved exactly in all the cases.

\begin{figure}[htbp]
\begin{center}
\begin{tabular}{cc}
\includegraphics[scale=0.35]{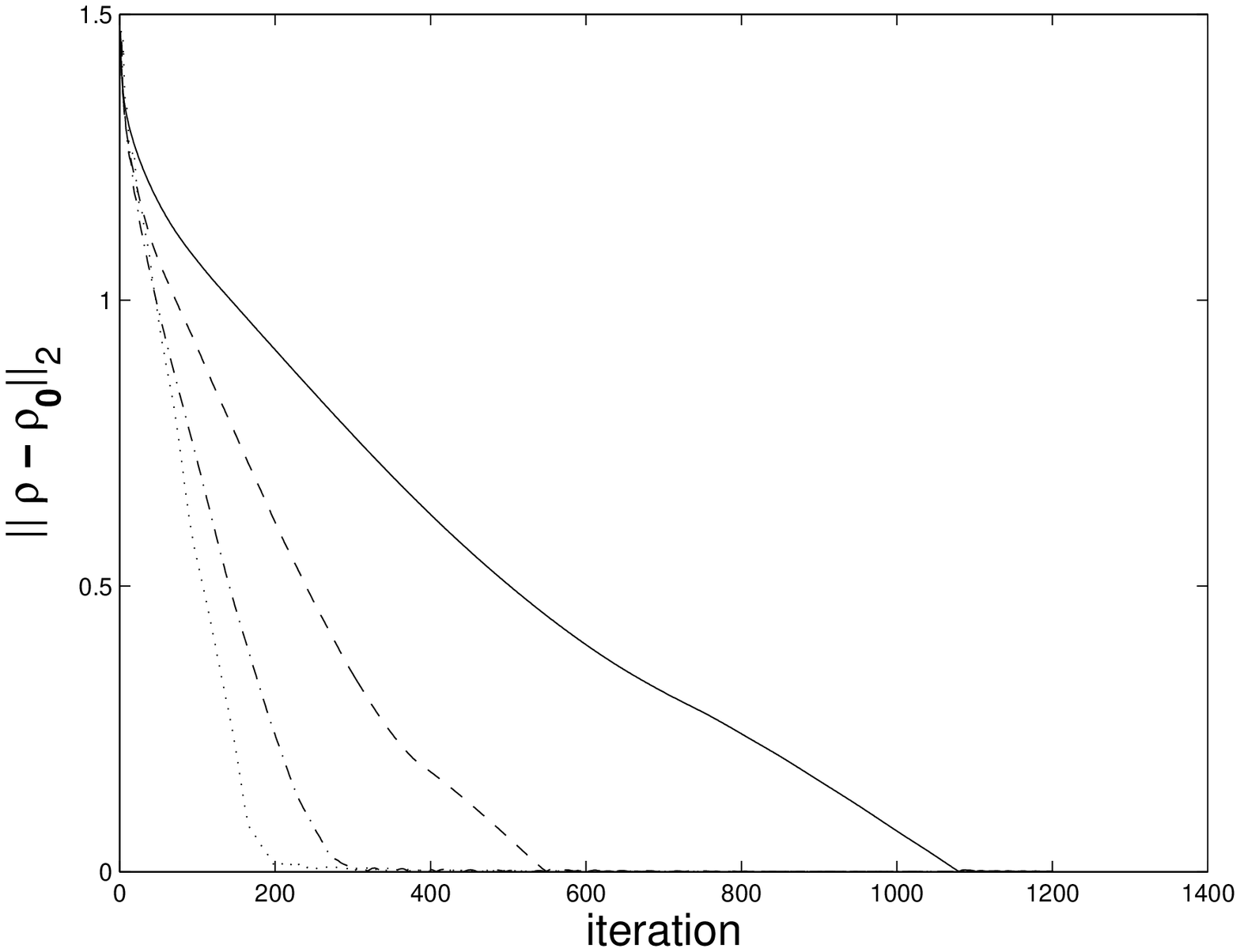}&
\includegraphics[scale=0.35]{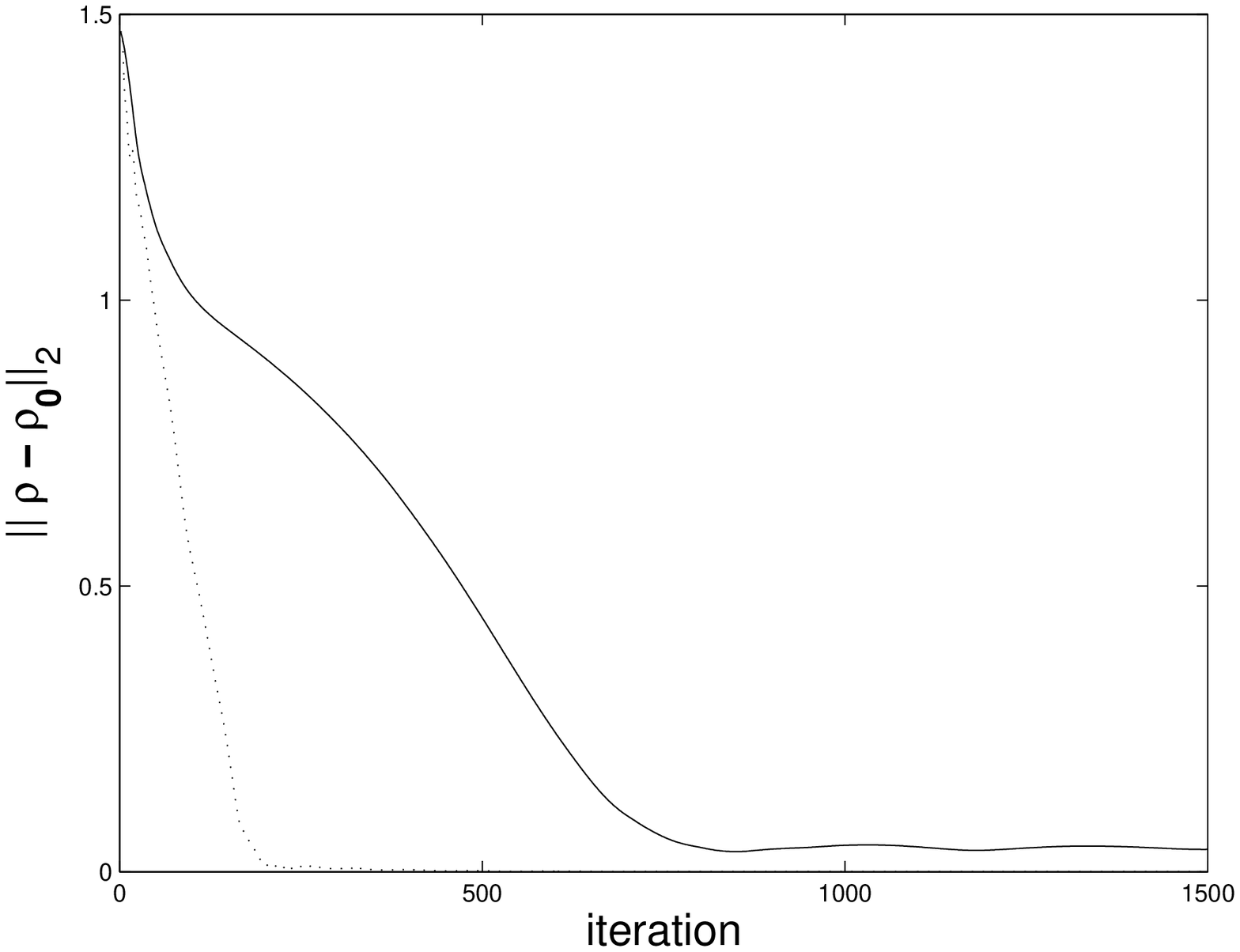}
\end{tabular}
\caption{Right: Plots of the convergence rate of the GeLMA
algorithm  for various values of $\tau=\alpha\,\|A_\omega^T{\bf
b}(\omega)\|_{l_\infty}$: $\alpha=2$ (solid line), $\alpha=5$
(dashed line), $\alpha=10$ (dot-dashed line), and $\alpha=20$
(dotted line). Left: Comparison of the converge rates of the GeLMA
algorithm with $\alpha=20$ (dotted line) and the FISTA method with
$\alpha=0.01$ (solid line). In these numerical experiments we have
used the four scatterers configuration shown in the top right
image of Fig. \ref{fig:imagesnonoise}. Noiseless data.}
\label{fig:comparison}
\end{center}
\end{figure}

An interesting feature of the GeLMA
algorithm \eqref{fd2bis2} is that it attains the exact
solution of the basis pursuit problem
for large values of the regularization parameter $\tau$.
This speeds up the convergence
rate.
Informally, this speed-up of convergence can
be seen from the coercivity estimate (\ref{jan1804}) and the
error estimate (\ref{jan1802}).
Note that for other popular gradient based algorithms, such as
ISTA or FISTA \cite{Beck09}, $\tau$ has to be smaller than
$\|A_\omega^T{\bf b}(\omega)\|_{l_\infty}$. Otherwise, they
converge to the (maximally sparse) zero solution $\bfrho = 0$. To
examine this property in more detail, we show in
Fig.~\ref{fig:comparison} (left panel) plots of the $\ell_2$
distance to the exact solution $\|\bfrho - \bfrho_0 \|$ as
function of the iteration number for various values of $\tau=
\alpha \|A_\omega^T{\bf b}(\omega)\|_{l_\infty}$: $\alpha=2$
(solid line), $\alpha=5$ (dashed line), $\alpha=10$ (dot-dashed
line), and $\alpha=20$ (dotted line). We observe that the larger
the value of $\tau$ is, the faster is the convergence rate.
Furthermore, for all the values of $\tau$ the algorithm achieves
the exact solution $\bfrho_0$.

In Fig. \ref{fig:comparison} (right panel) we compare the
convergence rates of the GeLMA algorithm and the FISTA algorithm
\begin{eqnarray}
\bfrho^{(k)} = \eta_{\tau \alpha_k} (\bfrho^{(k)} - \alpha_k \nabla f({\bf \xi}^{(k)}))\, , \label{fista1}\\
 \alpha_{k+1} = \frac{1 + \sqrt{1 + 4 \alpha^2_k}}{2} \, , \label{fista2}\\
 {\bf \xi}^{(k + 1)} = \bfrho^{(k)} + \frac{\alpha_{k} - 1}{\alpha_{k+1}} (\bfrho^{(k)} - \bfrho^{(k-1)})\, ,
\label{fista3}
\end{eqnarray}
for $\tau= 0.01 \|A_\omega^T{\bf b}(\omega)\|_{l_\infty}$. We
choose a small value of $\tau$ because we are considering
noisefree data in these examples. In
\eqref{fista1}-\eqref{fista3}, $\bfrho_{1}$ and ${\bf \xi}_2 =
\bfrho_{1}$ are given, and $\alpha_1< 2/L$. We observe that the
convergence rate of the FISTA algorithm (solid line) for $\tau=
0.01 \|A_\omega^T{\bf b}(\omega)\|_{l_\infty}$ is much slower than
the convergence rate of the GeLMA algorithm for $\tau= 20
\|A_\omega^T{\bf b}(\omega)\|_{l_\infty}$. Even more, the FISTA
algorithm  with $\tau= 0.01 \|A_\omega^T{\bf
b}(\omega)\|_{l_\infty}$ does not obtain the exact solution. To
achieve the exact solution, we would have to let $\tau \rightarrow
0$.

\begin{figure}[htbp]
\begin{center}
\begin{tabular}{ccc}
\includegraphics[scale=0.25]{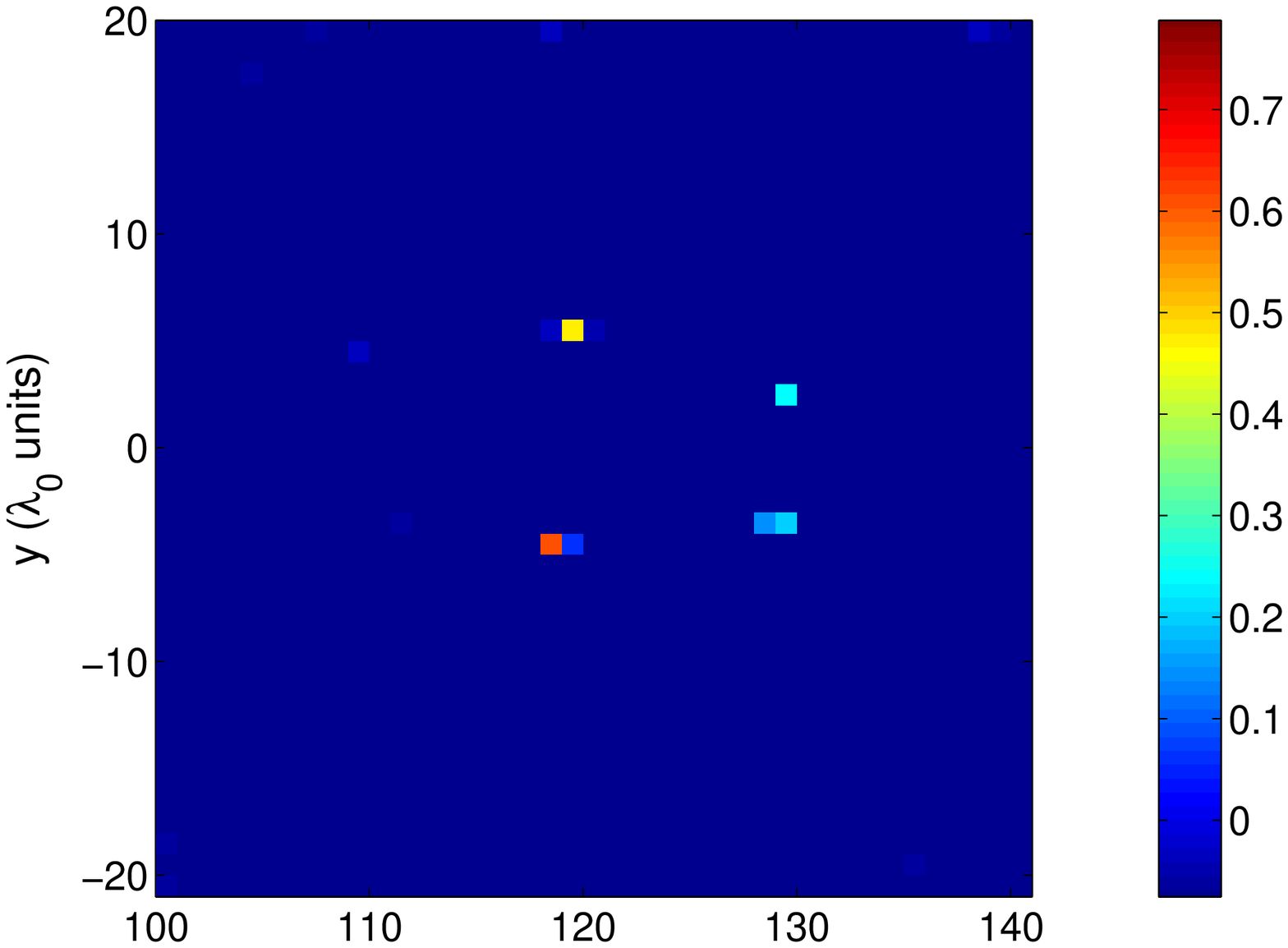}&
\includegraphics[scale=0.25]{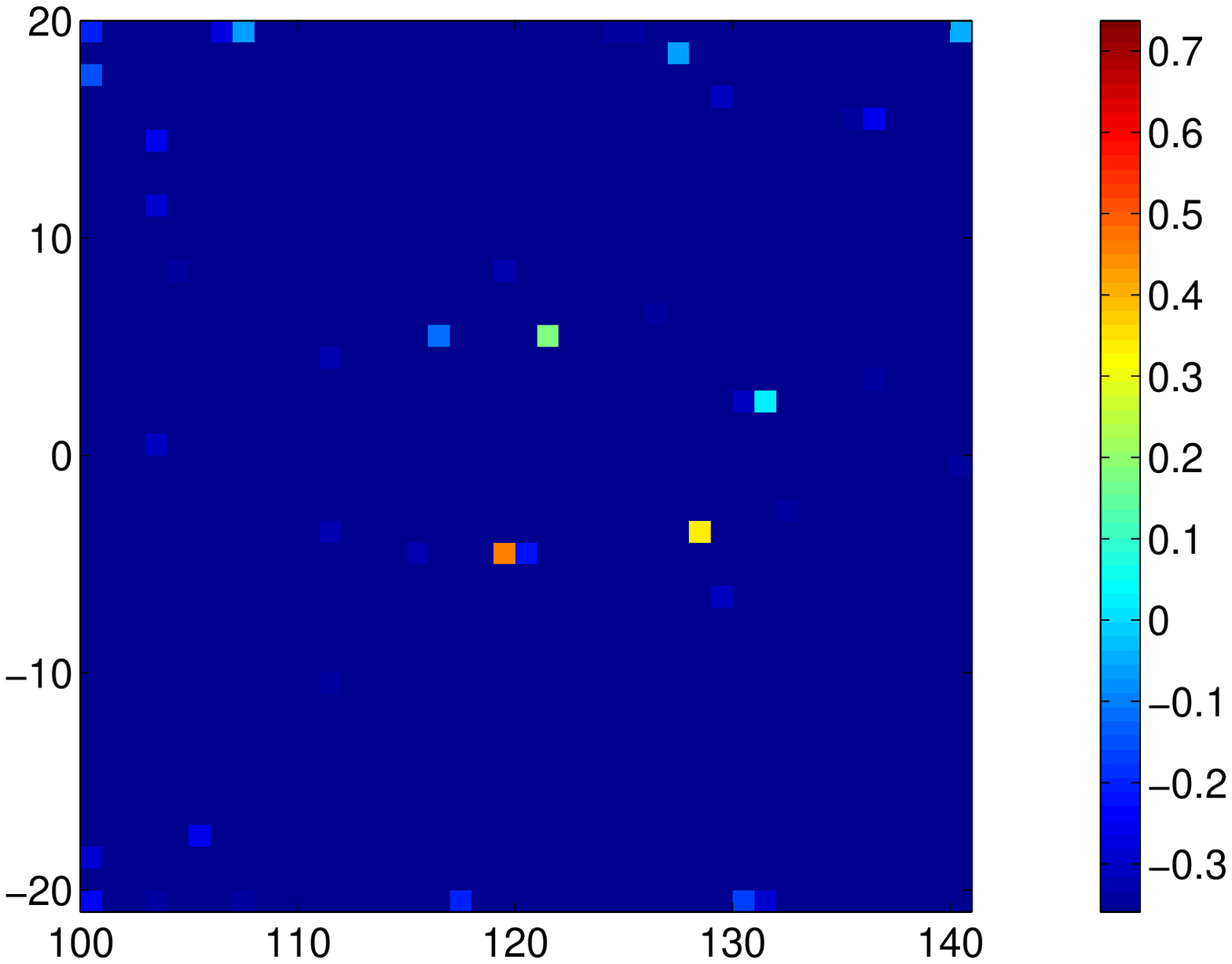}&
\includegraphics[scale=0.25]{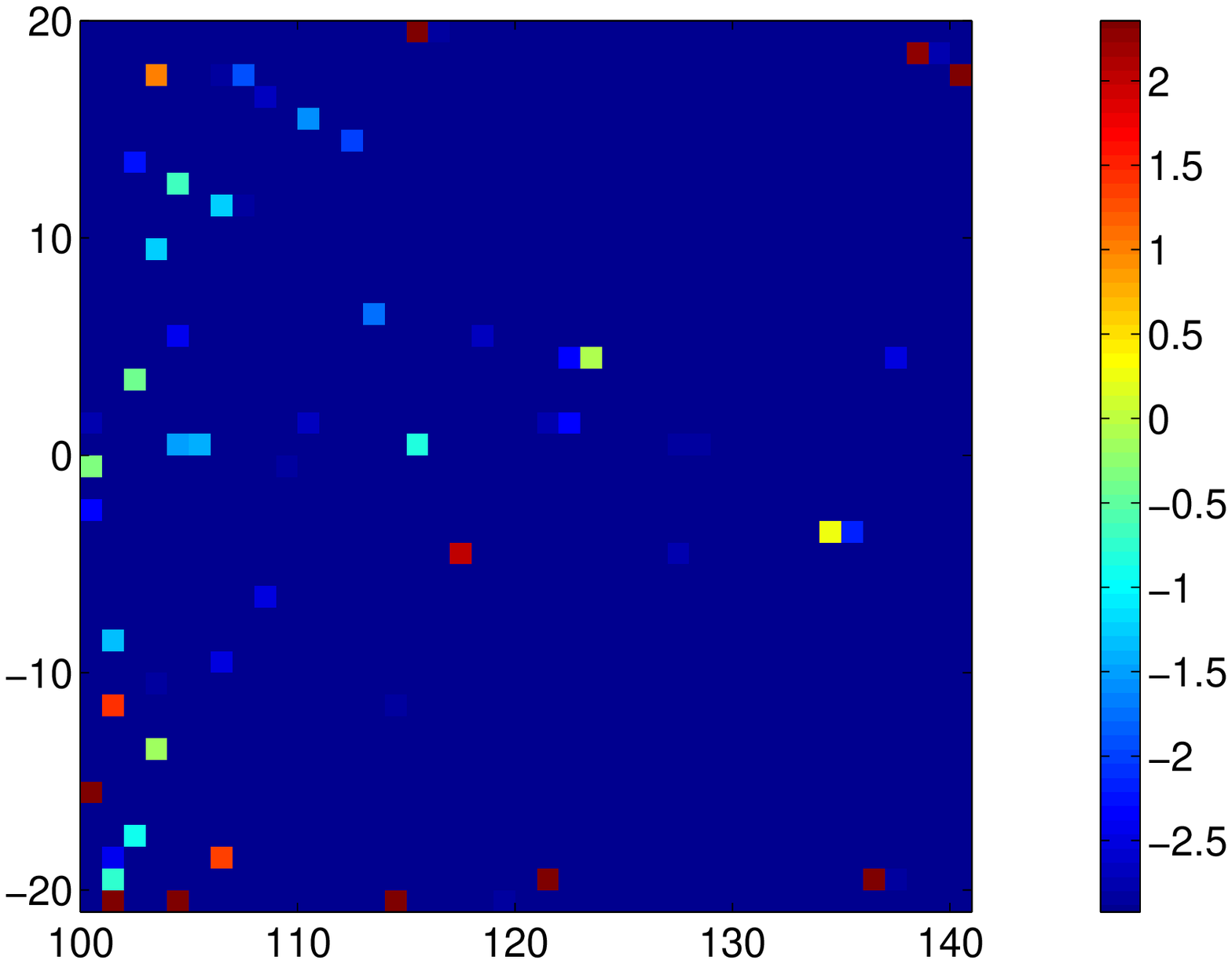}\\
\includegraphics[scale=0.25]{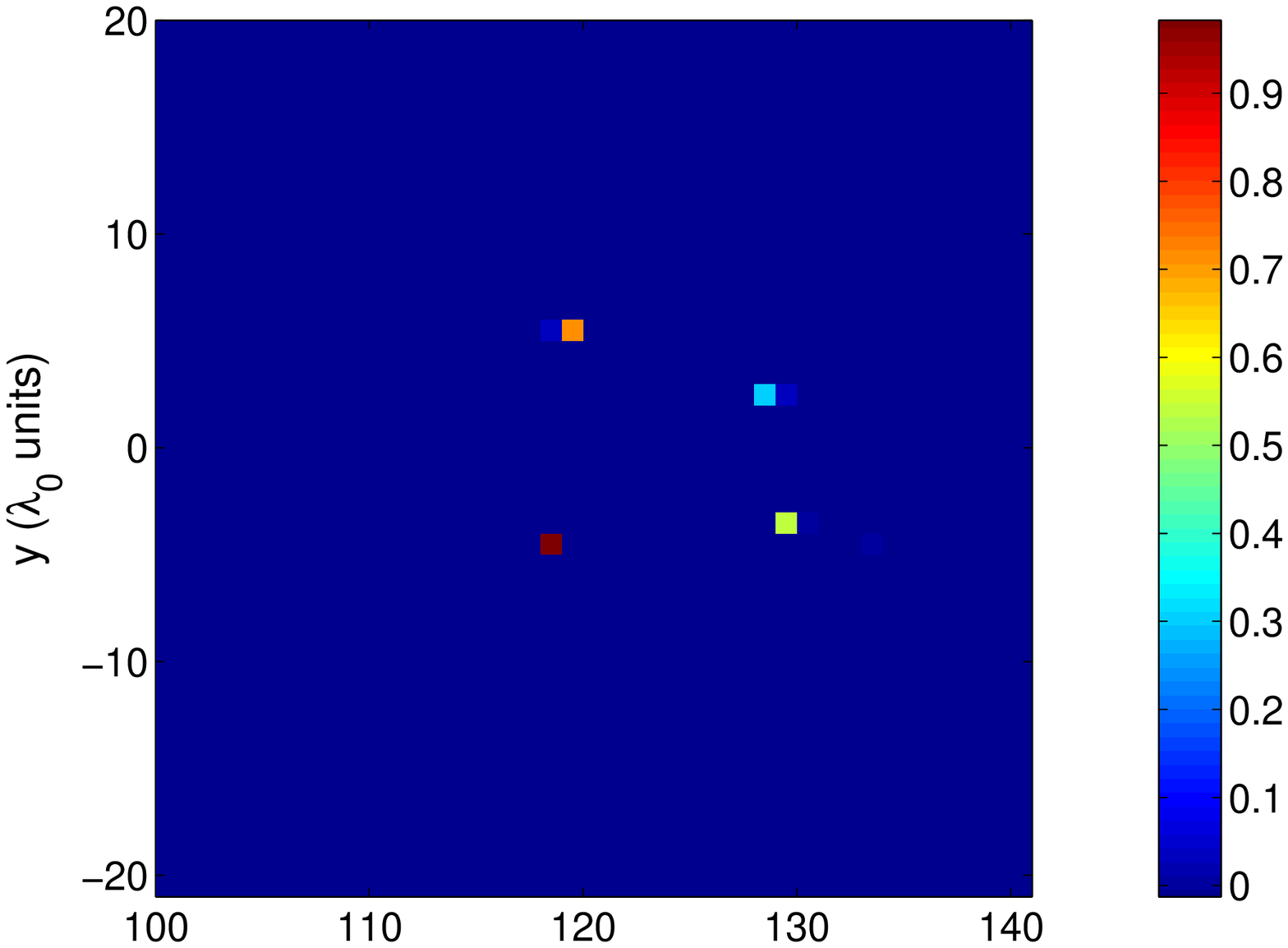}&
\includegraphics[scale=0.25]{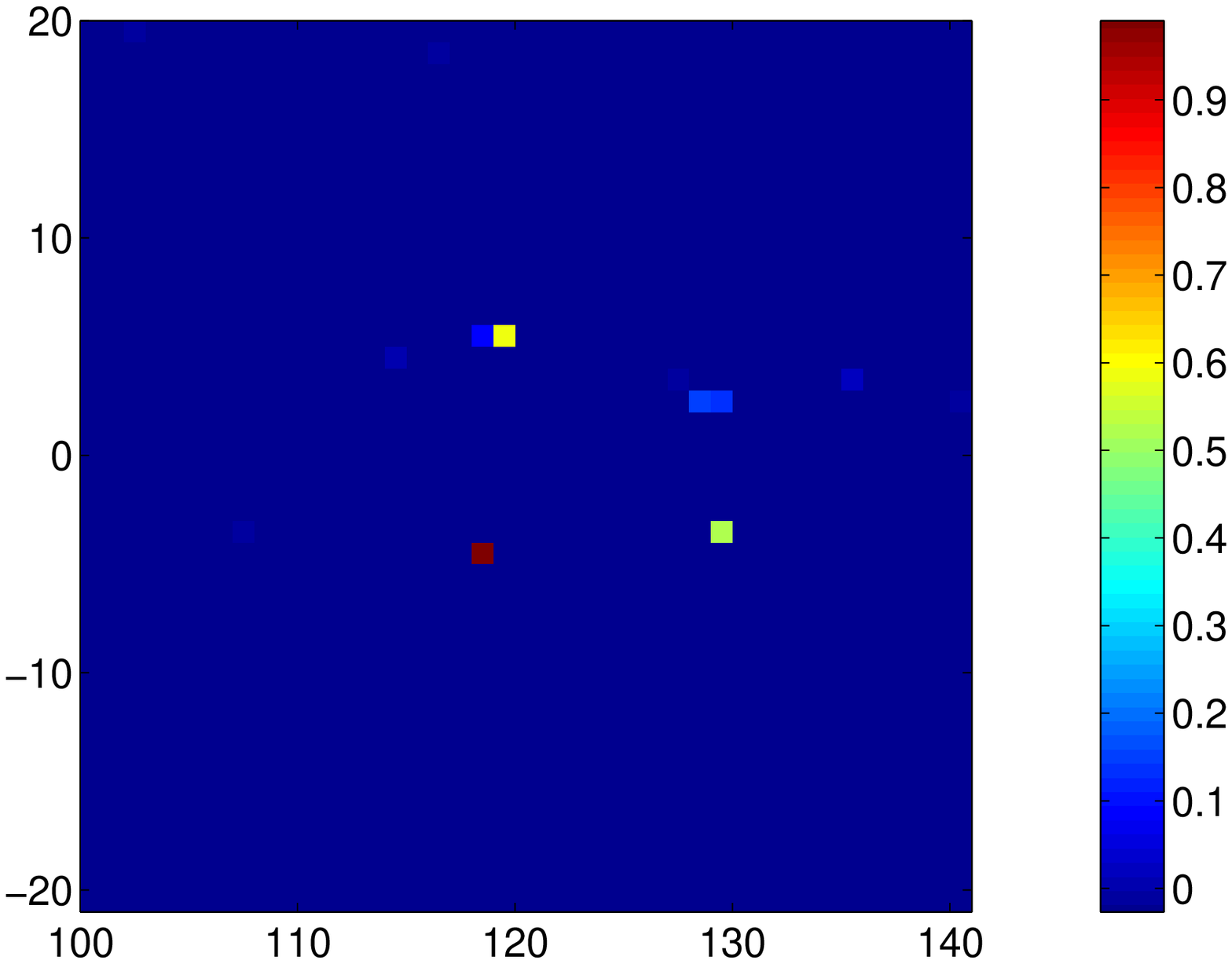}&
\includegraphics[scale=0.25]{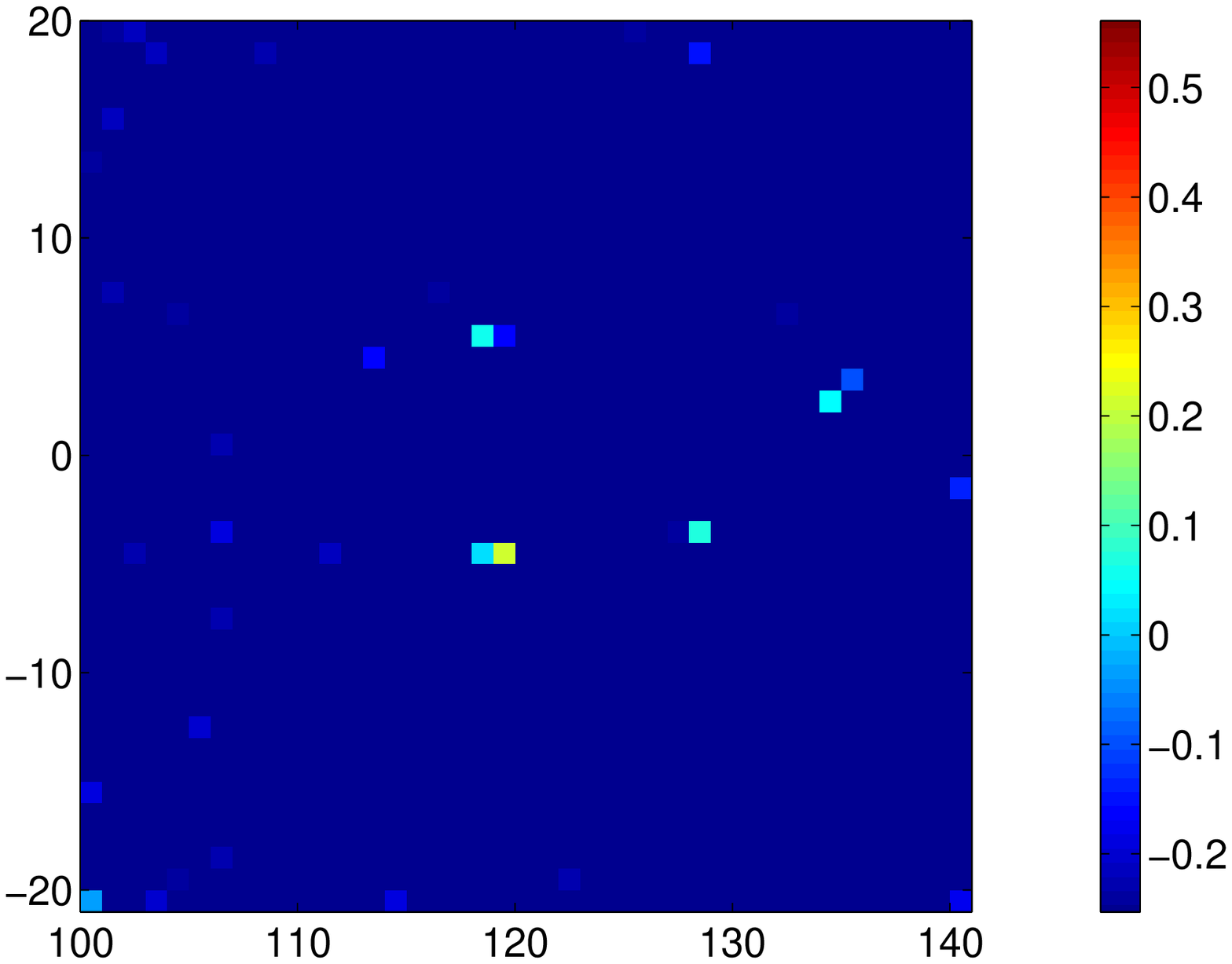} \\
\includegraphics[scale=0.25]{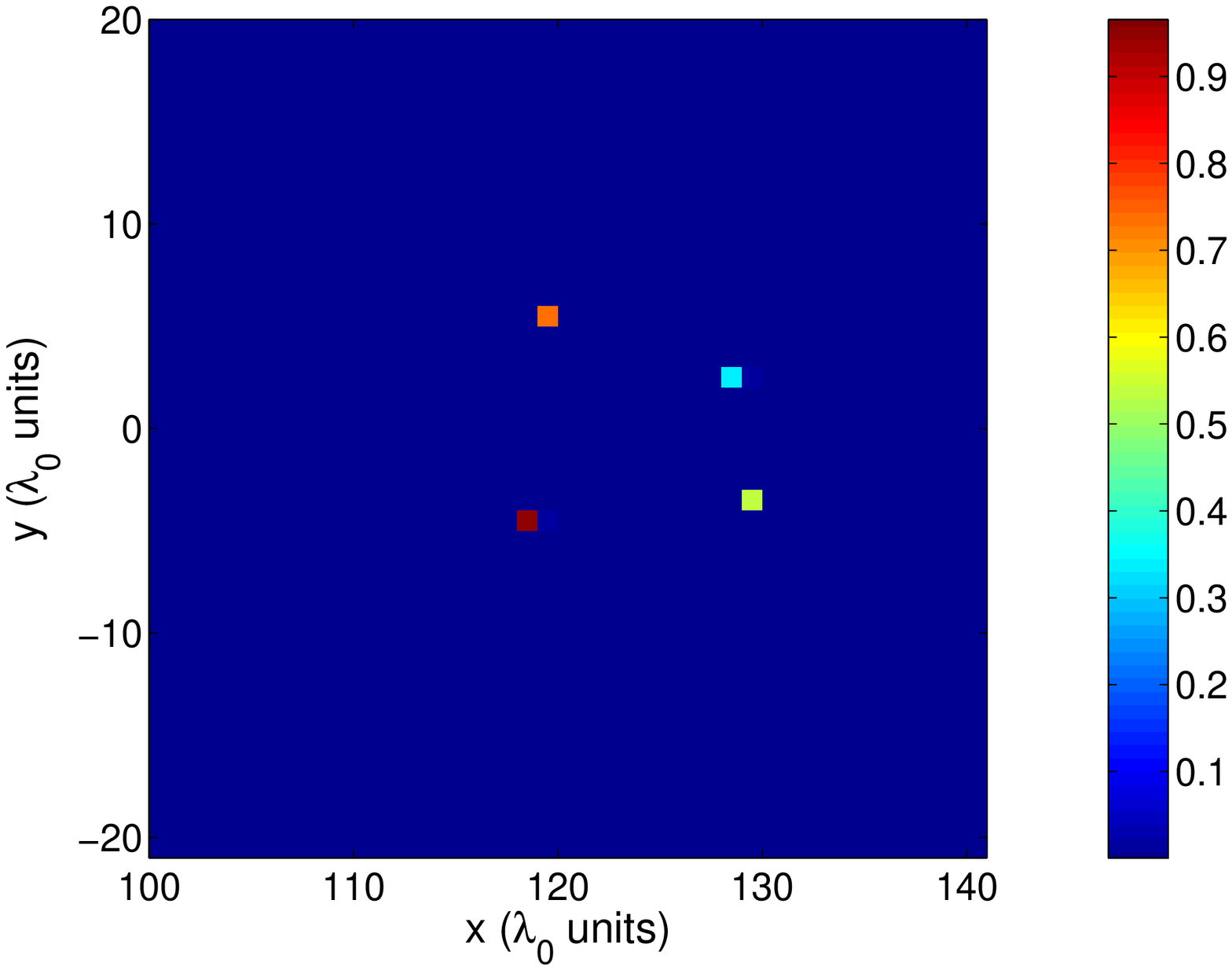}&
\includegraphics[scale=0.25]{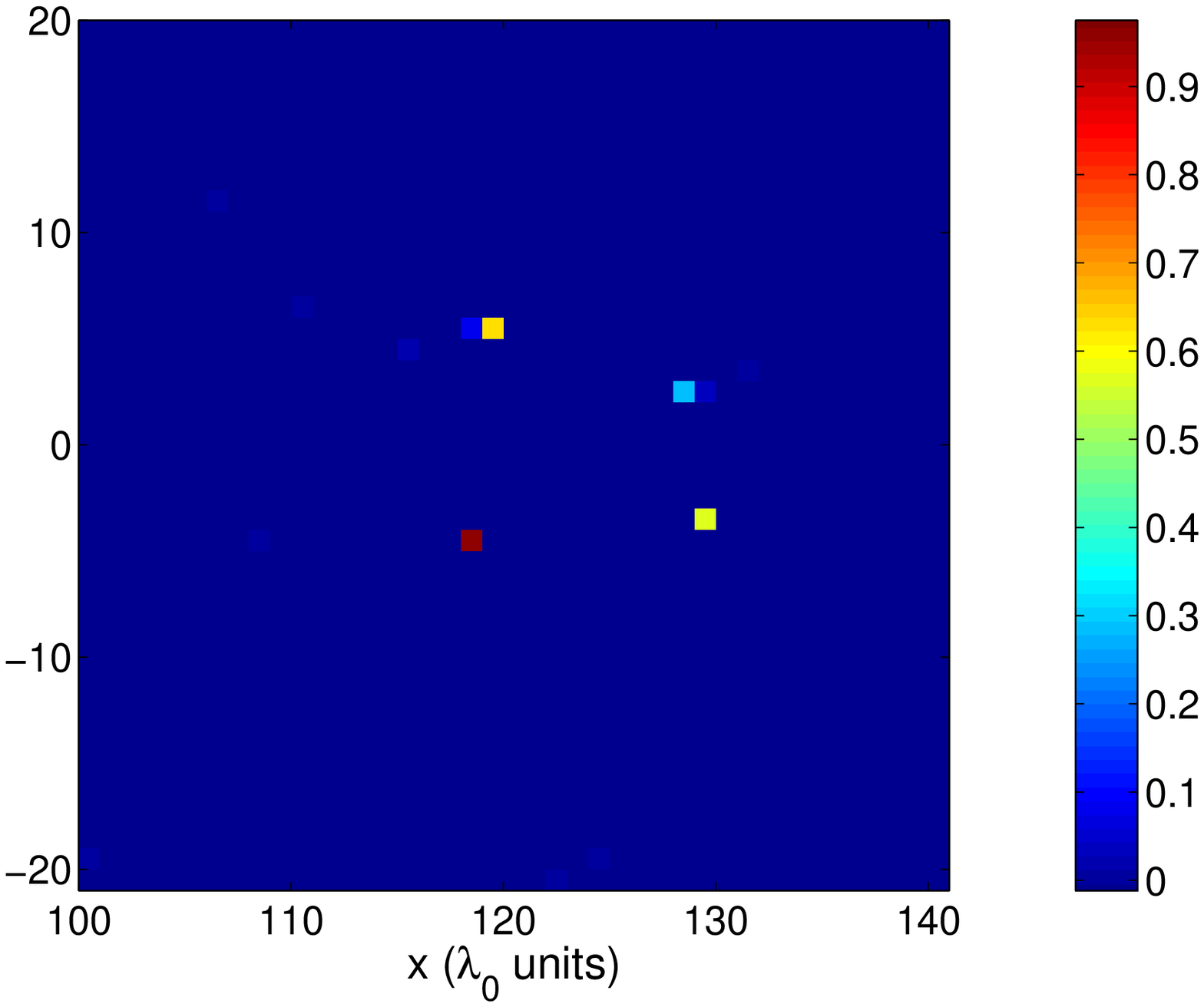}&
\includegraphics[scale=0.25]{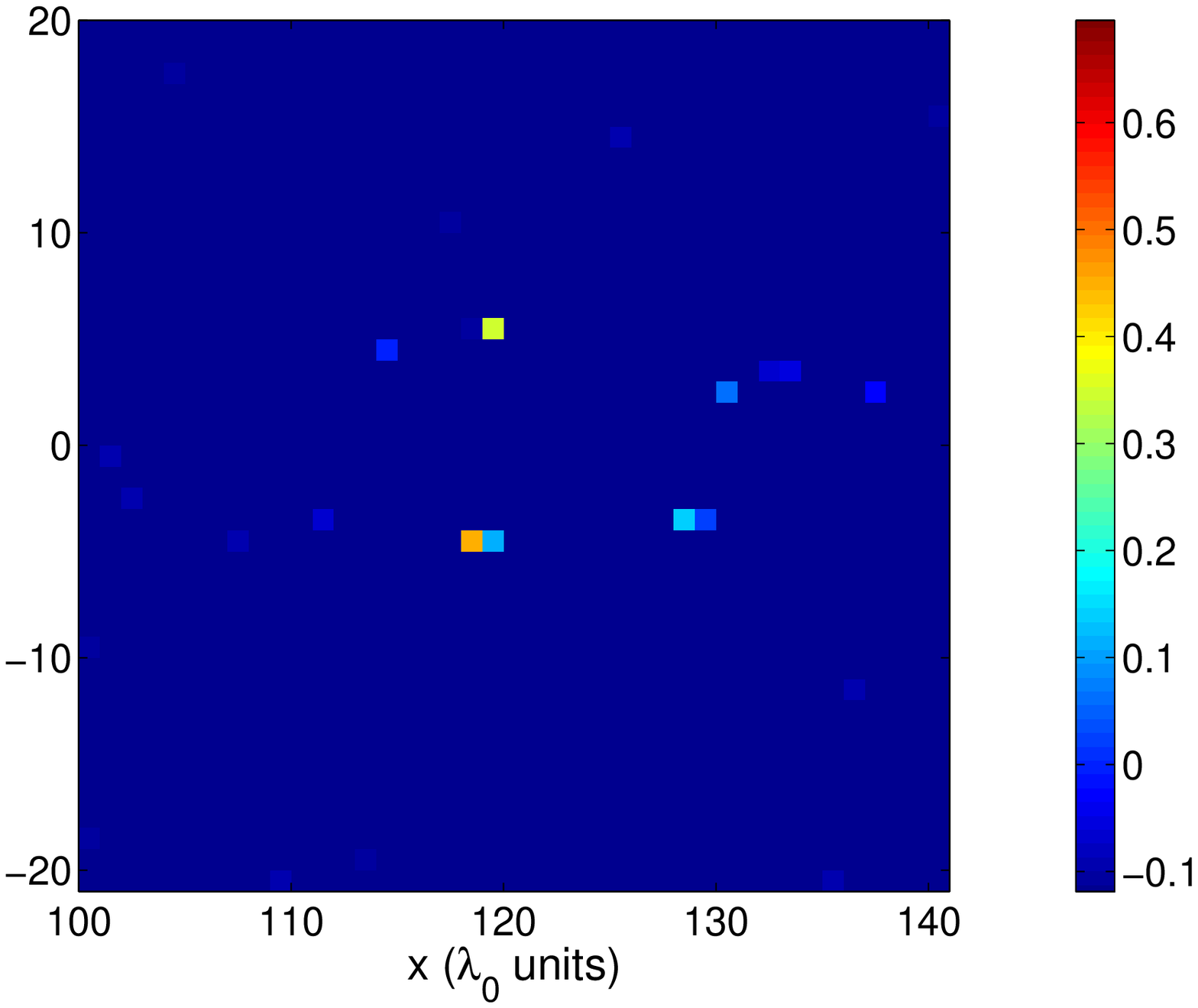}
\end{tabular}
\caption{Impact of the regularization parameter
$\tau=\alpha\,\|A_\omega^T{\bf b}(\omega)\|_{l_\infty}$ on the
reconstructions for different amounts of noise in the data. Top
row: Recovered images with $\alpha=2$ and $5\%$ noise (left),
$10\%$ noise (middle) and $30\%$ noise (right). Middle row: same
as top row but for $\alpha=20$. Bottom row: same as top row but
for $\alpha=200$. } \label{fig:imagesnoise}
\end{center}
\end{figure}

Next, we examine the performance of the GeLMA algorithm under
noise contaminated data $\vect b(\omega) + \vect e(\omega)$. The
noise vector $\vect e(\omega)$ is generated by independent
Gaussian random variables with zero mean and standard deviation
$\beta\,\|\vect b(\omega)\|/\sqrt{N}$. Here, $\beta$ is a
parameter that measures the noise strength. In Fig.
\ref{fig:imagesnoise}, we show the results for $\beta=0.05$ (left
column), $\beta=0.1$ (middle column), and $\beta=0.3$ (right
column). For a fixed step size $\Delta t$, the regularization
parameter $\tau=\alpha\,\|A_\omega^T{\bf b}(\omega)\|_{l_\infty}$
controls the sparsity of the solution. Hence, one expects the
algorithm to be more stable with respect to additive noise when
$\tau$ is large. We plot in Fig. \ref{fig:imagesnoise} the
recovered images using different values of $\tau$: $\alpha=2$ (top
row), $\alpha=20$ (middle row) and $\alpha=100$ (bottom row). We
observe in the top row that the location of the scatterers is
recovered exactly when there is $5\%$ noise in the data (left
plot). The recovered reflectivities are also quite close to the
real ones. However, when the noise increases to $10\%$ (middle
plot) one scatterer is missing in the recovered image that also
shows some ghost scatterers. As expected, the image gets worse
when the noise is $30\%$, as can be seen in the right plot. The
results are much better when we increase the value of $\alpha$ to
$20$ (middle row). With $5\%$ noise in the data (left plot) both
the location and reflectivities of the scatterers are very close
to the real ones.  Even with $10\%$ noise in the data (middle
plot) we can determine the location of the four scatterers.
However, with $30\%$ noise we miss the forth scatterer. Finally,
we show in the bottom row the recovered images using $\alpha=200$.
For $5\%$ and $10\%$ noise (left and middle images, respectively),
the location of the scatterers is exact. Furthermore, the
recovered reflectivities are very sharp. However, we still miss
the location of one scatterer when there is $30\%$ noise in the
data, as can be seen in the right image of the bottom row of this
figure. We plan to investigate in detail the robustness of the
algorithm with respect to noise in a future publication.

\section{Proof of Theorems~\ref{thm-jan10},~\ref{thm-sub} and~\ref{thm-ode}}
\label{sec:proofs1}

Theorems~\ref{thm-jan10} and \ref{thm-sub} are easy consequences of Theorem~\ref{thm-ode} and its proof.

\subsection{Outline of the proof of Theorem~\ref{thm-ode}}

Let $\bar x$ be the unique minimizer of (\ref{l1})
We write $x_\eps=\bar{x} + q_\eps$ and obtain
\begin{eqnarray} \label{dqbis}
\frac{dq_\eps}{dt} = -\tau G_{\eps}\left(\bar{x} + q_\eps \right)
+A^* \left( z_\eps -Aq_\eps \right),~~~~
\frac{dz_\eps}{dt} =  - Aq_\eps\,\label{s12} .
\end{eqnarray}
Our goal is now to show that $q_\eps(t)\to 0$ as $t\to+\infty$.
If we take the time-derivative of the first equation in (\ref{dqbis}), 
and use the second equation, we obtain:
\begin{eqnarray}\label{s21bis}
\ddot{q}_\eps + A^*A (\dot{q}_\eps+q_\eps) = - \tau g^{\eps} \left(\bar{x} + q_\eps \right)\dot{q}_\eps.
\end{eqnarray}
Here $g^\eps(x)$ is a diagonal matrix with the entries on the main diagonal given by
\begin{eqnarray}\label{nov806}
g_{ii}^\eps(x)=\begin{cases}
&0,~\hbox{ if } |x_i|>\eps,\\
&1/\eps,~\hbox{ if } |x_i|<\eps.\\
\end{cases}
\end{eqnarray}
Note that (\ref{s21bis}) is simply an equation for an oscillator with friction, and
a forcing term in the right side. As the matrix $A^*A$ is singular, 
the oscillator is degenerate.
Therefore, it is reasonable to expect that the   friction term
$A^*A\dot q_\eps$ in (\ref{s21bis}) by itself
would ensure that $Aq_\eps(t)\to 0$ as $t\to+\infty$, provided that the forcing does
not interfere. However, the friction alone can not send $q_\eps(t)$ to
zero since it is degenerate. 
Moreover, in showing that $q_\eps(t)$ becomes small as $t\to+\infty$,
one has to use the fact that $\bar x$ is the minimizer of (\ref{l1})
and not just any solution of $Ax=y$.  The strategy of the proof is (i)
to establish uniform bounds on $q_\eps(t)$ and $z_\eps(t)$, and (ii) show that
any limit point of $q_\eps(t)$ as $t\to+\infty$ is close to zero. 

The a priori bounds are obtained in several steps. We first describe the required 
intermediate
lemmas, and present their proofs later.
The first step in the proof is the following lemma that provides a
Lyapunov function for (\ref{dqbis}) and establishes a bound on
$\|Aq_\eps(t)\|$.
\begin{lemma}\label{lemnov82}
There exists a constant $C_0>0$ that is independent of $\eps$ (and depends only
on the initial data $x_0$) so that
\begin{equation}\label{nov808}
\|\dot q_\eps(t)\|^2+\|Aq_\eps(t)\|^2+\int_0^\infty \|A\dot q_\eps(s)\|^2ds<C_0,
\end{equation}
for all $\eps<\eps_0$ andf all $t>0$.
\end{lemma}
The bound on $\|A\dot q_\eps\|$ in Lemma~\ref{lemnov82} leads to a uniform bound on $z_\eps(t)$.
\begin{lemma}\label{lemnov83}
There exists a constant $C>0$ that is independent of $\eps>0$ so that
$\|z_\eps(t)\|\le C$ for all $t>0$.
\end{lemma}
The next step is to show that $Aq_\eps(t)$ is small for large times.
Since $\dot z_\eps=Aq_\eps$, it follows from Lemma~\ref{lemnov83} that
\[
\int_{t_1}^{t_2}Aq_\eps(s)ds
\]
is uniformly bounded for all $t_{1,2}>0$. Together with the
integral bound on $A\dot q_\eps(t)$ in Lemma~\ref{lemnov82}, this shows that
$Aq_\eps(t)$ becomes small at some "not too large" time.
\begin{lemma}\label{lemdec6-4}
There exists two constants $C_{1,2}>0$ that are independent of
$\eps\in(0,\eps_0)$ so that for any $k\in\Nm$ there exists a
time $t_k<C_1k^{3}$
such that for all $t\in(t_k,t_k+C_2k)$ we have
$\|Aq_\eps(t)\|\le C_1/k$ for all $\eps<\eps_0$.
\end{lemma}
Next, using the bounds in Lemmas~\ref{lemnov82} and~\ref{lemnov83}, as
well as the precise form of the forcing term in (\ref{s21bis}), we obtain
a uniform bound for $\|q_\eps(t)\|$:
\begin{lemma}\label{propdec7-6}
There exists a constant $C>0$ so that   we have
\begin{equation}\label{dec704}
 \|\bar x+q_\eps(t)\| \le {C},
\end{equation}
for all $t>0$ and all $\eps>0$.
\end{lemma}
The bound on $\|q_\eps(t)\|$ allows us to strengthen Lemma~\ref{lemdec6-4} to
include a bound on $\dot q_\eps(t)$ "at some times" as well.
\begin{lemma}\label{cor-jan3}
There exists a constant $C>0$ that is independent of $\eps\in(0,\eps_0)$ so that for any $k\in\Nm$ there exists a
time $s_k<Ck^{3}$
such that
$\|Aq_\eps(s_k)\|^2+\|\dot q_\eps(s_k)\|^2\le C/k$ for all $\eps<\eps_0$.
\end{lemma}
The Lyapunov function in Lemma~\ref{lemnov82} and Lemma~\ref{cor-jan3}  together
imply that $\dot q_\eps(t)$ and $\dot z_\eps(t)$ are not only
"small sometimes" but rather tend to zero
as $t\to+\infty$ 
\begin{corollary}\label{cor-jan3-2}
There exists a constant $C>0$ that is independent of $\eps\in(0,\eps_0)$ so that for any $n\in\Nm$ there exists a
time $s_n=s_n(\eps)<Cn^{3}$
such that
$\|Aq_\eps(s)\|^2+\|\dot q_\eps(s)\|^2\le C/n$ for all $\eps<\eps_0$ and all $s>s_n$.
\end{corollary}
Corollary~\ref{cor-jan3-2} shows that the right side of the ODE system (\ref{dqbis})
is small as $t\to+\infty$. The final step in the proof is to show that this implies
that $q_\eps(t)$ is small, and it is here that the condition that $\bar x$
is the minimizer of (\ref{l1}) comes into play.

\subsection{The end of the proof of Theorem~\ref{thm-ode} }

It follows from Corollary~\ref{cor-jan3-2} that for any $\delta_0>0$ there exist $T=T(\delta_0)$, and  $\eps_0=\eps_0(\delta_0)$
\begin{equation}\label{jan401}
 \| A^* z_\eps(t) - \tau G_\eps(\bar{x}+q_\eps(t))  \| \leq \delta_0, \| A q_\eps(t) \| \leq \delta_0
\end{equation}
for all $\eps \leq \eps_0$ and $t>T$.
The first inequality in~\eqref{jan401} implies
\[
\left| (A^* z_\eps(t) - \tau G_\eps(\bar{x}+q_\eps(t) ) \cdot (\bar{x}+q_\eps(t)) \right| \leq \delta_0 \| \bar{x} +q_\eps(t) \|.
\]
Using the second inequality from~\eqref{jan401} in
\[
\left| A^* z_\eps(t) \cdot (\bar{x}+q_\eps(t)) - A^* z_\eps(t) \cdot \bar{x} \right| \leq \delta_0 \| z_\eps (t)  \|,
\]
and denoting\footnote{The quantity $\aleph_\eps(x)$ plays essentially the same role as
$\| x \|_{l^1_\eps}$ defined in~\eqref{nov826}. They are, however, quantitatively  slightly different.}
\[
\aleph_\eps(x)= \sum_i x_i G_\eps(x_i),
\]
we obtain
\[
\left| \tau \aleph(\bar{x} +q_\eps(t)) - A^* z_\eps(t) \cdot \bar{x} \right| \leq \delta_0 \left( \| z_\eps (t) \|  +\| \bar{x} +q_\eps(t) \|  \right).
\]
It also follows from the first inequality in~\eqref{jan401}  that
\[
\| A^* z_\eps(t) \|_{l_\infty}\leq \tau+\delta_0,
\]
and thus
\[
\left| (A^* z_\eps(t) \cdot \bar{x} )\right| \leq (\tau+\delta_0) \| \bar{x} \|_{l_1}.
\]
As a consequence,
\[
\aleph(\bar{x} +q_\eps(t)) - \| \bar{x} \|_{l_1}
\leq \frac{\delta_0}{\tau} \left( \| \bar{x} \|_{l_1}  +\| z_\eps (t) \|  +\| \bar{x} +q_\eps(t) \|  \right).
\]
and therefore
\begin{equation}\label{jan1802}
 \| \bar{x} +q_\eps(t) \|_{l_1}  - \| \bar{x} \|_{l_1}  \leq \frac{\delta_0}{\tau} \left( \| \bar{x} \|_{l_1}  +\| z_\eps (t) \|  +\| \bar{x} +q_\eps(t) \|  \right)+\eps_0 n.
\end{equation}
Here $n$ is the dimension of $q_\eps$. As $\bar x$ is the unique minimizer,
for any $\delta$ we can choose $\alpha$ and $\delta_0$ sufficiently small~so that estimates
\[
\| \bar{x} +q_\eps(t) \|_{l_1}  - \| \bar{x} \|_{l_1}  \leq \alpha,~ \| A q_\eps(t) \| \leq \delta_0
\]
imply that $\|q_\eps\|<\delta$.
Hence it remains to use uniform boundedness of  $\bar{x} +q_\eps(t)$ and $z_\eps(t)$ and choose $\delta_0$ and $\eps_0$ so that
\[
\frac{\delta_0}{\tau} \left( \| \bar{x} \|_{l_1}  +\| z_\eps (t) \|  +\| \bar{x} +q_\eps(t) \| \right) +\eps_0 n \leq \alpha.
\]
This finishes the proof of Theorem~\ref{thm-ode} except for the proof of 
Lemmas~\ref{lemnov82}-\ref{cor-jan3} and Corollary~\ref{cor-jan3-2}.~$\Box$

 \subsection{Proof of Theorem~\ref{thm-jan10}}

Fix $T_\delta$ such that $|q_\eps(t)|<\delta$ for all $T>T_\delta$.
We know from the Arzela-Ascoli theorem that $q_\eps(t)\to q(t)$ and $z_\eps\to z(t)$ uniformly on $[0,T_\delta]$,
after extracting a subsequence,
and the functions $q(t)$ and $z(t)$ are Lipschitz on $[0,T_\delta]$, with the Lipschitz constant independent of $\delta>0$.
The second equation in (\ref{dqbis}), and the dominated convergence theorem imply that
\begin{equation}\label{jan1006}
z(t)=-\int_0^tAq(s)ds,
\end{equation}
whence
\begin{equation}\label{jan1008}
\dot z=-Aq,~~z(0)=0.
\end{equation}

The family $f_\eps(t)=G_\eps(\bar x+q_\eps(t))$ is uniformly bounded in $L^2[0,T_\delta]$. Hence, after possibly extracting a
subsequence, it converges weakly in $L^2[0,T_\delta]$ to a limit $f(s)$. The (vector-valued) function $f(s)$ satisfies the following
properties: (i) $-1\le f_j(t)\le 1$, for all $0\le t\le T_\delta$, $1\le j\le N$, and (ii) if $q_j(t)\neq -\bar x_j$ then $f_j(t)=\hbox{sgn}(\bar x_j+q_j)$.
It follows that for any $0\le t_1<t_2\le T_\delta$ we have
\begin{equation}\label{jan1010}
q(t_2)-q(t_1)=-\tau \int_{t_1}^{t_2}f(s)ds+\int_{t_1}^{t_2} A^*(z(s)-Aq(s))ds.
\end{equation}
The aforementioned properties  of $f(t)$ imply that $x(t)=\bar x+q(t)$ is a strong solution of (\ref{jan602}).
Uniqueness of the strong solution~\cite{CL} implies that the whole family $x_\eps(t)=\bar x+q_\eps(t)$, $z_\eps(t)$ converges
to the solution of (\ref{jan602}). The conclusion of Theorem~\ref{thm-jan10} now follows from Theorem~\ref{thm-ode}.~$\Box$

\subsection{Proof of Theorem~\ref{thm-sub}}

Theorem~\ref{thm-ode} implies that as $\eps \to 0$ and $t \to \infty$, along a subsequence, we have
$\bar z_{\eps_k}\to \bar z$ and $\bar q_{\eps_k}\to \bar 0$.
Then the first estimate in~\eqref{jan401} implies that $\bar\lambda=\bar z/\tau$ satisfies
\begin{eqnarray}\label{nov922}
&&[A^*\bar \lambda]_j= \hbox{sgn}\bar x_j ,\hbox{ if $\bar x_j\neq 0$}.
\\
&&-1\le [A^*\bar \lambda]_j\le 1,\hbox{ if $\bar x_j= 0$}.\nonumber
\end{eqnarray}
This completes the proof of Theorem~\ref{thm-sub}.~$\Box$

\section{Proofs of auxiliary lemms for the proof of Theorem~\ref{thm-ode}}

\subsection{Proof of Lemma~\ref{lemnov82}}
Multiplying~\eqref{s21bis} by $\dot{q}_\eps(t)$, gives
\begin{equation}\label{lyapunov0bis}
\frac{1}{2}\frac{d}{dt} \left( \|\dot{q}_{\eps}(t)\|^2 +
\|A q_{\eps}(t)\|^2 \right) = -\|A \dot{q}_{\eps}(t)\|^2 -
\tau \langle g^\eps(\bar{x}+q_\eps) \dot q_\eps,\dot{q}_\eps\rangle.
\end{equation}
Let
\begin{equation}\label{ntbis}
N^\eps_t =\int_0^t   \langle g^\eps(\bar{x}+q_\eps(s)) \dot q_\eps(s),\dot{q}_\eps(s)\rangle ds=
\sum_{j=1}^n\int_0^t g_{jj}^\eps(\bar{x}_j+q_{\eps,j}(s)) |\dot q_{\eps,j}(s)|^2ds\ge 0,
\end{equation}
then integrating~\eqref{lyapunov0bis} in time we get
\begin{equation}\label{lyapunovbis}
\frac{1}{2}\left( |\dot{q}_{\eps}(0)|^2 +
\|A q_{\eps}(0)\|^2\right) = \frac{1}{2}\left( \|\dot{q}_{\eps}(T)\|^2 +
\|A q_{\eps}(T)\|^2\right) + \tau N^\eps_T + \int_{0}^T\|A \dot{q}_{\eps}\|^2 dt,
\end{equation}
and (\ref{nov808}) follows.
Note that $\|\dot q_\eps(0)\|$ is uniformly bounded
in $\eps>0$ since the function $G_\eps(s)$ takes values
in the interval~$[-1,1]$.~$\Box$

\subsection {Proof of Lemma~\ref{lemnov83}}
Differentiating the second equation in (\ref{dqbis}) we obtain
\begin{equation}\label{s1_bis}
\ddot{z}_\eps +AA^*(\dot{z}_\eps+z_\eps) =\tau A G_\eps  \left( \bar{x}+q_\eps(t) \right).
\end{equation}
Let us multiply this equation by $e^t$ and integrate, to obtain
\begin{equation}\label{nov824}
\int_{0}^{t} e^s \ddot{z}_\eps(s) ds + e^tAA^*z_\eps (t) = \tau A \int_{0}^{t}e^{s} G_\eps  \left( \bar{x}+q_\eps(s) \right) ds,
\end{equation}
since $z(0)=0$.
We estimate, using (\ref{nov808}):
\[
\left\| \int_{0}^{t} e^t \ddot{z}_\eps(s) ds \right\| \leq \left( \int_{0}^{t} e^{2t}  ds  \int_{0}^{t}  \|\ddot{z}_\eps(s)\|^2 ds\right)^{1/2}
\leq
C e^t\left(\int_0^t\|A\dot q_\eps(s)\|^2ds\right)^{1/2}\le Ce^t.
\]
As $|G_{\eps,j}|\le 1$ for all $1\le j\le n$, we also have
\[
\left\| A \int_{0}^{t} e^{s} G_\eps  \left( \bar{x}+q_\eps(s) \right) ds  \right\| \leq C e^t.
\]
Since the matrix $AA^*$ is invertible, we obtain from (\ref{nov824}) that
$\|z_\eps(t)\|  \leq C$.~$\Box$

\subsection{Proof of Lemma~\ref{lemdec6-4}} 

Let us set $y_\eps(t)=Aq_\eps(t)$.
As $z_\eps(t)$ is uniformly bounded,  there exists a constant $C>0$ that is independent of $\eps$ so that
\begin{equation}\label{dec604}
\int_{t_1}^{t_2}y_\eps(s)ds<C,
\end{equation}
for all $0<t_1<t_2$.
If we take an integer $n=Ck$, we have
\begin{equation}\label{dec616}
 \left\|\frac{1}{n}\int_{t}^{t+n}y_\eps(s)ds\right\|<\frac{1}{2k},
\end{equation}
for all $t>0$,
and
\begin{equation}\label{dec6-18}
\|y_\eps(t)\|\le \left\|\frac{1}{n}\int_{t}^{t+n}y_\eps(s)ds\right\|+\frac{1}{n}\int_t^{t+n}\!\!\!\sqrt{s-t}
\left(\int_t^s\|\dot y_\eps(\xi)\|^2d\xi\right)^{1/2}\!\!\! ds
\le \frac{1}{2k}+ \sqrt{n}\left(\int_t^{t+n}\|\dot y_\eps(s)\|^2ds\right)^{1/2}\!\!\!.
\end{equation}
Lemma~\ref{lemnov82} implies that given $n$ there exist at
most $Ck^2{n}=Ck^{3}$ integers $l$ such that
\[
\int_l^{l+2n}\|\dot y_\eps(s)\|^2ds>\frac{1}{4k^2{n }}.
\]
It follows that there exists $k_0<Ck^{3}$ such that
\[
\int_{k_0}^{k_0+2n}\|\dot y_\eps(s)\|^2ds<\frac{1}{4k^2{n}}.
\]
Then, for all $t\in(k_0,k_0+n)$ we have
\[
\int_{t}^{t+n}\|\dot y_\eps(s)\|^2ds<\frac{1}{4k^2{n}},
\]
whence
\begin{equation}\label{dec6-20}
\|y_\eps(t)\|\le \farc{C}{k},
\end{equation}
for all $t\in(k_0,k_0+n)$.~$\Box$

\subsection{Proof of Lemma~\ref{propdec7-6}} 

Let us recall (\ref{s21bis}):
\begin{eqnarray}\label{dec710}
\ddot{q}_\eps + A^*A (\dot{q}_\eps+q_\eps) = - \tau g^{\eps} \left(\bar{x} + q_\eps \right)\dot{q}_\eps,
\end{eqnarray}
Multiply this equation by $q_\eps$ and integrate:
\begin{eqnarray}\label{dec712}
&&\langle q_\eps(t),\dot q_\eps(t)\rangle-\langle q_\eps(0),\dot q_\eps(0)\rangle+\frac{1}{2}\|Aq_\eps(t)\|^2-\farc{1}{2}\|Aq_\eps(0)\|^2
+\int_0^t\|Aq_\eps(s)\|^2ds\nonumber\\
&&=\int_0^t\|\dot q_\eps(s)\|^2ds
-\tau \int_0^t  \langle g^{\eps} \left(\bar{x} + q_\eps \right)\dot{q}_\eps,q_\eps\rangle ds
\end{eqnarray}

Next, set
\[
v_\eps(t)=-\int_0^tq_\eps(s)ds,
\]
so that  $z_\eps(t)=Av_\eps(t)$, and $v_\eps(0)=0$. We rewrite~\eqref{dqbis} as
\begin{eqnarray} \label{dqlbis}
\frac{dq_\eps}{dt}= -\tau G_{\eps}\left(\bar{x} + q_\eps \right) +A^*A \left( v_\eps - q_\eps \right),~~~
\frac{dv_\eps}{dt} =  - q_\eps  \label{s12l}.
\end{eqnarray}
Consider the function
\begin{equation}\label{jan306}
Q(t)= \frac{1}{2} \| A(v_\eps(t)-q_\eps(t))\|^2+\tau  \| \bar{x} +q_\eps(t)\|_{l^1_\eps}.
\end{equation}
Then we have
\begin{eqnarray}\label{jan304}
&&\frac{d Q}{dt} = \tau\langle  G_{\eps}\left(\bar{x} + q_\eps \right),\dot q_\eps\rangle
-  \langle A^*A \left( v_\eps - q_\eps \right),\dot{q}_\eps \rangle +\langle A^*A \left( v_\eps - q_\eps \right),\dot{v}_\eps \rangle
\\
&&= -  \|\dot{q}_\eps \|^2 +
\frac{1}{2}\frac{d}{dt} \|Av_\eps \|^2 -\langle A \dot{v}_\eps, A q_\eps\rangle
=-  \|\dot{q}_\eps \|^2 +\frac{1}{2}\frac{d}{dt} \|z_\eps \|^2 +\|\dot{z}_\eps \|^2.\nonumber
\end{eqnarray}
As $z_\eps(0)=0$, it follows that
\begin{eqnarray*}
 &&\| A(v_\eps(t)-q_\eps(t))\|^2+\tau  \| \bar{x} +q_\eps(t)\|_{l^1_\eps} -  \| A(v_\eps(0)-q_\eps(0))\|^2-
 \tau  \| \bar{x} +q_\eps(0)\|_{l^1_\eps}
 \\
 &&
   = \frac{\|z_\eps(t)\|^2}{2}  - \int_0^t  \| \dot{q}_\eps(s) \|^2ds + \int_0^t \| \dot{z}_\eps(s)\|^2 ds.
\end{eqnarray*}
This can be re-written as
\begin{equation}\label{dec708bis}
 \| z_\eps(t) -Aq_\eps(t) \|^2+\tau  \| \bar{x} +q_\eps(t)\|_{l^1_\eps} + \int_0^t  \|\dot{q}_\eps(s)\|^2ds =
 \frac{\|z_\eps(t)\|^2}{2} + \int_0^t \|A{q}_\eps(s)\|^2 ds+C_0.
\end{equation}
Adding (\ref{dec712}) and (\ref{dec708bis}) gives:
\begin{eqnarray}\label{dec714}
&&\langle q_\eps(t),\dot q_\eps(t)\rangle+\frac{1}{2}\|Aq_\eps(t)\|^2
+ \| z_\eps(t) -Aq_\eps(t)\|^2+\tau  \| \bar{x} +q_\eps(t)\|_{l^1_\eps} +  \nonumber\\
&&= -\tau \int_0^t  \langle g^{\eps} \left(\bar{x} + q_\eps \right)\dot{q}_\eps,q_\eps\rangle ds+
\frac{\|z_\eps(t)\|^2}{2} + C_0',
\end{eqnarray}
with the constant $C_0'$ that only depends on the initial data. Lemmas~\ref{lemnov82} and~\ref{lemnov83} imply then
\begin{eqnarray}\label{dec716}
\langle q_\eps(t),\dot q_\eps(t)\rangle +\tau  || \bar{x} +q(t)||_{l^1_\eps} =
-\tau \int_0^t  \langle g^{\eps} \left(\bar{x} + q_\eps \right)\dot{q}_\eps,q_\eps\rangle ds+
r(t),
\end{eqnarray}
with a uniformly bounded function $r(t)$: $|r(t)|\le C$. We claim that there exists $C>0$ that is independent of $\eps$ and $t$ so that
\begin{equation}\label{dec718}
\left| \int_0^t  \langle g^{\eps} \left(\bar{x} + q_\eps \right)\dot{q}_\eps,q_\eps\rangle ds\right|\le C.
\end{equation}
Indeed, let us fix some $1\le j\le n$ and look at
\begin{equation}\label{dec720}
I=\int_0^t    g_{jj}^{\eps} \left(\bar{x} + q_\eps (s)\right){q}_{\eps,j}(s)\dot q_{\eps,j}(s)  ds=
\farc{1}{\eps}\sum_{k=1}^Q\int_{s_k}^{s_k'} {q}_{\eps,j}(s)\dot q_{\eps,j}(s)  ds=
\farc{1}{2\eps}\sum_{k=1}^Q  (\|{q}_{\eps,j}(s_k')\|^2-\| q_{\eps,j}(s_k)\|^2).
\end{equation}
Here $(s_k,s_k')$, $k=1,\dots,Q$, are the time intervals that $q_j(s)$ spends in the interval $(-\bar x_j-\eps,-\bar x_j+\eps)$,
and $q_j(s_k)=\bar x_j\pm \eps$, depending on whether $q_j$ enters this interval from above or below, and similarly for $q_\eps(s_k')$.
It is easy to see that $q_j(s_k')=q_j(s_{k+1})$, whence (\ref{dec720}) is a telescoping sum, giving
\[
I= \farc{1}{2\eps}   (\|{q}_{\eps,j}(s_Q')\|^2- \| q_{\eps,j}(s_1)\|^2).
\]
As both terms in the right side above can take only the values $-\bar x_j\pm \eps$, we conclude that $|I|\le C$, so that (\ref{dec718})
holds. Now, (\ref{dec716})  becomes
\begin{equation}\label{dec724}
\langle q_\eps(t),\dot q_\eps(t)\rangle +\tau  \| \bar{x} +q_\eps(t)\|_{l^1_\eps} \le C.
\end{equation}
As $\| x\| \leq \| x \|_{l^1}$, using the traingle inequality, we obtain the following inequality for $m_\eps(t)=\|q_\eps(t)\|$:
\[
m_\eps(t)\dot m_\eps(t)+C_1m_\eps(t)\le C_2.
\]
Now, the comparison principle implies that $m_\eps(t)\le C'$ for all $t>0$, and the proof of Lemma~\ref{propdec7-6} is complete.~$\Box$

\subsection{Proof of Lemma~\ref{cor-jan3}} 

Let us choose $t_k$ and $t_{k'}$ as in the proof of Lemma~\ref{lemdec6-4}.
The estimate for $\|Aq_\eps(t_k)\|$ is exactly as in that Lemma.
%
Next, dividing (\ref{dec708bis}) by $C_2k=t_k'-t_k$ we get, due to the boundedness of $z_\eps(t)$ and $q_\eps(t)$:
\begin{equation}\label{jan308}
\farc{1}{t_k'-t_k}\int_{t_k}^{t_k'}\|\dot q_\eps(s)\|^2ds\le
\frac{C}{k}+\frac{C}{k}\int _{t_k}^{t_k'}\|Aq_\eps(s)\|^2ds\le\frac{C}{k}+\farc{C}{k^2}.
\end{equation}
It follows that
there exists  a time $s_k\in(t_k,t_k')$ such that
$\|\dot{q}_\eps(s_k)\|\le C/\sqrt{k}$.~$\Box$
 
 \subsection{Proof of Corollary~\ref{cor-jan3-2}} 

This follows immediately from Lemma~\ref{cor-jan3} and (\ref{lyapunov0bis}), as the latter implies that
\begin{equation}\label{lyapunovbis17}
\|\dot{q}_{\eps}(t)\|^2 + \|A q_{\eps}(t)\|^2 \le  \|\dot{q}_{\eps}(s_n)\|^2 + \|A q_{\eps}(s_n)\|^2\le\farc{C}{n},
\end{equation}
for all $t>s_n$.~$\Box$

\commentout{
\subsubsection{Improved estimates in the limit $\eps \to 0$}
Note that we used equivalence of $l^2$ and $l^1$ norms after~\eqref{dec724} in order to conclude~\eqref{dec704}. Thus the constant in~\eqref{dec704} may a priori depend on the dimension. Here we show that we can only use the existence of
$\bar{\lambda}$ from~\eqref{nov922} in order to
obtain~\eqref{dec704} where the constant is independent of dimension. We will also denote
\[
G(\bar{x} +q(t) ):=\lim_{\eps \to 0} G_\eps(\bar{x} +q_\eps(t)),~G(\bar{x}):=A^*\lambda,~w(t)=z(t)-\tau\bar{\lambda}.
\]
Then~\eqref{jan602} becomes
\[
\dot{q}=-\tau (G(\bar{x}+q) - G(\bar{x}))+ A^*(w- Aq),~\dot{w}=-Aq.
\]
Multiplying the first equation by $q$, the second by $w$ and integrating we obtain
\[
\|q(t)\|^2 +\|w(t)\|^2 +\tau \int_0^t \|Aq(s)\|^2 ds + \tau   \int_0^t q(t) \cdot \left( G(\bar{x}+q(t)) - G(\bar{x}) \right)ds = \|q(0)\|^2 +\|w(0)\|^2,
\]
where we note that $q(t) \cdot \left( G(\bar{x}+q(t)) - G(\bar{x}) \right) \geq 0$ by monotonicity of $G$. In particular we conclude
\[
 \int_0^t \|Aq(s)\|^2 ds \leq C.
\]
Combining the last inequality with~\eqref{dec708bis} we conclude
\[
\| \bar{x} +q(t)\|_{l^1} \leq C,~\int_0^t  \|\dot{q}\|^2ds \leq C,
\]
where $C$ is independent of dimension.

}

\section{The proof of Theorem~\ref{thm-minimexact}}
\label{sec:proofs2}

We will use Theorem~\ref{thm-sub} in order to prove Theorem~\ref{thm-minimexact}.
The role of the vector $z$ that satisfies the sub-differential condition can be seen from the following lemma.

\begin{lemma}\label{lem-subdiff}
Suppose the sub-differential condition does not hold for a particular ${z}$. Then for this $z$ we have a strict inequality
\begin{equation}\label{nov1404}
h({z})= \min_{x} F(x,{z})  < \tau || \bar{x}||_{l_1}.
\end{equation}
\end{lemma}
{\bf Proof.}
Assume that $z$ does not satisfy the sub-differential condition, that is, either
\begin{itemize}
\item[(i)] $\left| \left[ A^*{z} \right]_i\right| > \tau$ for some $i$, or
 \item[(ii)] $\left| \left[ A^*{z} \right]_i\right| \leq \tau$, but $\left[ A^*{z} \right]_i   \neq  \tau \hbox{sign}(\bar{x}_i)$ for
 some $i$ such that $\bar{x}_i \neq 0$.
\end{itemize}
We will
show that
\begin{equation}\label{smallerbis}
F(\bar{x} +q,{z}) < F(\bar{x},{z}) = \tau ||\bar{x} ||_{l_1},
\end{equation}
for some (sufficiently small) $q$, which implies (\ref{nov1404}). We will now construct  $q$ explicitly.

Assume first that (i) holds:  $ |\left[ A^*{z} \right]_i|> \tau$ for some $i$. Then, set $r=\left[ A^*\tilde{z} \right]_i$, and
choose $q$ so that
\[
q_k=\begin{cases}
&\varepsilon   \hbox{ sign}\left(r \right), ~\hbox{ if } k=i,\\
&0, \hbox{ otherwise.}
\end{cases}
\]
We have
\begin{eqnarray*}
&&
\!\!\!\!\!\!\!\!\!\!\!\!\!\!\!
F(\bar{x} +q,{z}) = \tau\|\bar x+q\|_{l_1}+\farc{1}{2}\|Aq\|^2-\langle A^* z,q\rangle=
 \tau || \bar{x} ||_{l_1}  + \tau ( |\bar{x}_i + \varepsilon  \hbox{sign}\left(r \right)|- |\bar{x}_i|  ) +   \frac{1}{2}\| Aq\|^2 -
 \varepsilon  |r|\\&&\le \tau || \bar{x} ||_{l_1}+\eps\tau-\eps| r|  +\frac{1}{2}\| Aq\|^2\le \eps \tau || \bar{x} ||_{l_1}+\eps\tau-\eps| r| +C\eps^2
<  \tau || \bar{x} ||_{l_1},
\end{eqnarray*}
provided that we choose $\eps$ sufficiently small.

Similarly,  if (ii) holds, pick some $i$ such that   $\bar{x}_i \neq 0$ but $\left[ A^*{z} \right]_i \neq   \tau \hbox{sign}(\bar{x}_i)$.
Assume first that $\left[ A^*{z} \right]_i =   r \, \hbox{sign}(\bar{x}_i)$ with $0 < \left| r \right| < \tau$.
Pick $\varepsilon\in(0, |  \bar{x}_i | )$
and choose $q$ with the components
\begin{equation}\label{lastqbis}
q_k=\begin{cases}
&-\varepsilon \,  \hbox{sign}(\bar{x}_i), ~\hbox{ if } k=i,\\
&0, \hbox{ otherwise.}
\end{cases}
\end{equation}
The computation is similar:
\begin{eqnarray}\label{nov1002}
&& \!\!\!\!\!\!\!\!\!\!\!\!\!\!\! F(\bar{x} +q,{z}) = \tau\|\bar
x+q\|_{l_1}+\farc{1}{2}\|Aq\|^2-\langle  A^*\tilde  z,q\rangle=
 \tau || \bar{x} ||_{l_1}  + \tau ( |\bar{x}_i -\varepsilon  \hbox{sign}\left(\bar x_i \right)|- |\bar{x}_i|  )
 \\&&
  +   \frac{1}{2}\| Aq\|^2 +
 \varepsilon  r\le \tau || \bar{x} ||_{l_1}-\eps\tau+\eps r  +\frac{1}{2}\| Aq\|^2\le \eps \tau || \bar{x} ||_{l_1}-\eps\tau+\eps r+C\eps^2
<  \tau || \bar{x} ||_{l_1}, \nonumber
\end{eqnarray}
provided that $\eps$ is sufficiently small.
The last case case to consider is when (ii) holds,
but
$\left[ A^*{z} \right]_i  = - \tau  \hbox{sign}(\bar{x}_i)$. We still choose $q$ as in~\eqref{lastqbis},
and the computation is identical to (\ref{nov1002}), with $r=-\tau$. This completes the proof of Lemma~\ref{lem-subdiff}.~$\Box$
%
%

{\bf Proof of Theorem~\ref{thm-minimexact}.}
We trivially have
\[
h(z)=\min_{x} F(x,z) \leq F(\bar{x},z) =\tau ||\bar{x}||_{l_1},
\]
for all $z$.
 Thus, the conclusion of Theorem~\ref{thm-minimexact} would follow
 if we  show that there exists $\bar{z}$ such that $h(\bar z)=\tau\|\bar x\|_{l_1}$.
 That is, we need to show that for any $q \neq 0$ and some $\bar{z}$, we have
\begin{equation}\label{deffbis}
F(\bar{x}+q, \bar{z}) =  \tau || \bar{x}+q ||_{l_1} + \frac{1}{2}
\|Aq\|^2 - \langle A^* \bar{z}, q \rangle >  F(\bar{x},\bar{z})
=\tau ||\bar{x}||_{l_1}.
\end{equation}
We claim that~\eqref{deffbis} is true for any $\bar z$ that satisfies the sub-differential condition (\ref{nov1006}) --
recall that Theorem~\ref{thm-sub}
implies that  such $\bar z$
exists.
Let $\bar{z}$  satisfy the sub-differential condition~\eqref{nov1006}:  
\begin{eqnarray}\label{mainbcbis}
&&[A^* \bar{z}]_i =  \tau \hbox{ sign} \bar{x}_i,  \hbox{ if $i\in S_1$,}\\
&&\left| \left[ A^* \bar{z} \right]_i \right| \leq \tau, \hbox{ if $i\in S_0$.}\label{mainbcbis2}
\end{eqnarray}
We denoted here by $S_1$ the set of indices $i$ such that $\bar x_i\neq 0$, and by $S_0$ the set of indices $i$ such that $\bar x_i= 0$.

The function $F(\bar x+q,\bar{z})$ is convex in $q$. Hence, it suffices to show that $q=0$ is a strict local minimum, that is,
show that (\ref{deffbis}) holds for $q$ small enough.
In particular, we may assume that
\begin{equation}\label{tQbis}
\hbox{ sign} \left( \bar{x}_i  + q_i \right) =\hbox{ sign} \left( \bar{x}_i  \right),~~\hbox{if $i\in S_1$.}
\end{equation}
Now,   we obtain from (\ref{mainbcbis2}):
\begin{equation}\label{zero1bis}
\tau |{q}_i| -  \left[A^* \bar{z} \right]_i  {q}_i \geq 0, ~~i\in S_0,
\end{equation}
while for $i\in S_1$, we use~\eqref{tQbis} and~\eqref{mainbcbis} to obtain
\begin{equation}
 \tau |\bar{x}_i+q_i | -  \left[A^* \bar{z} \right]_i q_i =  \tau \, (\hbox{sgn}{\bar{x_i}}) (\bar{x}_i+q_i) -
 \tau  \,( \hbox{sgn}{\bar{x}_i}) q_i = \tau ( \hbox{sgn} {\bar{x_i}}) \bar{x}_i =\tau\, |\bar{x}_i|,~~i\in S_1.
 \label{zero2bis}
\end{equation}
We deduce from (\ref{zero1bis})-(\ref{zero2bis}) that
\begin{eqnarray}\label{nov1408}
&&\!\!\!
F(\bar x+q,\bar z)=  \tau || \bar{x}+q ||_{l_1} + \frac{1}{2} | Aq|^2 - \langle A^* \bar{z}, q \rangle  =
\sum_{i\in S_1}(\tau | \bar{x}_i+q_i |-[A^*z]_iq_i)+\sum_{i\in S_0}(\tau | q_i |-[A^*z]_iq_i)\nonumber\\
&&+ \frac{1}{2} | Aq|^2
\geq  \sum_{i\in S_1}\tau|\bar x_i|+ \frac{1}{2} | Aq|^2 =\tau ||\bar{x}||_{l_1}+ \frac{1}{2} | Aq|^2.
\end{eqnarray}
Therefore, we have $F(\bar x+q,\bar z)>\tau\|x\|_{l_1}$ unless
$Aq=0$. However, if $Aq=0$, then
\[
 F(\bar{x}+q,z) = \tau ||\bar{x}+q ||_{l_1} > \tau ||\bar{x} ||_{l_1},
\]
because $\bar{x}$ is the unique minimizer of~\eqref{l1}. 
Therefore, \eqref{deffbis} holds for all $q$.~$\Box$

\section{Conclusions}
\label{sec:conclusions}
We have shown using ordinary differential
equation methods that the relaxed $l_1$ minimization algorithm 
for problems with underdetermined linear constraints 
converges independently of the regularization parameter.
In the examples in array imaging the observed convergence
rates are faster than the theory implies, which means that more
analysis is needed. The algorithm is robust to noise although
we have not shown this theoretically. Finally, as the convergence rates
are independent of dimension, generalization to the infinite-dimensional case is
straightforward. 



{\small

 \end{document}